\documentclass[draft]{amsart}
\usepackage{amsmath}
\usepackage{amssymb}
\usepackage{graphicx}
\newtheorem{theorem}{Theorem}[section]
\newtheorem{lemma}[theorem]{Lemma}

\newtheorem{corollary}[theorem]{Corollary}
\theoremstyle{definition}

\newtheorem{example}[theorem]{Example}

\newtheorem{remark}[theorem]{Remark}

\numberwithin{equation}{section}
\begin{document}

\title[Quantization for Markov-type measures with complete overlaps]{Asymptotics of the quantization errors for some Markov-type measures with complete overlaps}
\author{Sanguo Zhu}

\address{School of Mathematics and Physics, Jiangsu University
of Technology\\ Changzhou 213001, China.}
\email{sgzhu@jsut.edu.cn}


\begin{abstract}
Let $\mathcal{G}$ be a directed graph with $2N$ vertices $1,2,\ldots, 2N$. Let $\mathcal{T}=(T_{i,j})_{(i,j)\in\mathcal{G}}$ be a family of contractive similitudes on $\mathbb{R}^q$. We
 assume that $T_{\widetilde{i},\widetilde{j}}=T_{i,j}$ for every $(\widetilde{i},\widetilde{j})$ in $\{(i,j),(i,j+N),(i+N,j),(i+N,j+N)\}\cap\mathcal{G}$. We denote by $K$ the Mauldin-Williams fractal determined by
 $\mathcal{T}$. Let $\chi$ be a positive probability vector and $P$  a $2N\times 2N$ row-stochastic matrix. We denote by $\nu$ the Markov-type measure associated with $\chi$ and $P$. Let $\mu$ be the image measure of $\nu$ under the natural projection, which is supported on $K$. We consider the following two cases: 1. $\mathcal{G}$ has two strongly connected components
 consisting of $N$ vertices; 2. $\mathcal{G}$ is strongly connected. With some assumptions for $\mathcal{G}$ and $\mathcal{T}$, for case 1, we determine the
 quantization dimension for $\mu$ in terms of the spectral radius of a related matrix; we prove that the lower quantization coefficient is always positive and establish a
 necessary and sufficient condition for the upper one to be finite. For case 2, we express the quantization dimension in terms of a pressure-like function and prove that
 the upper and lower quantization coefficient are always positive and finite.
\end{abstract}

\maketitle

\section{Introduction}

Let $\nu$ be a Borel probability measure on $\mathbb{R}^q$. The quantization problem for $\nu$ is concerned with the approximation of $\nu$ by discrete measures of finite support in
$L_r$-metrics. This problem has a deep background in information theory (cf. \cite{BW:82,GN:98}).  Rigorous mathematical foundations of quantization theory can be found in Graf and Luschgy's book \cite{GL:00}.

In the past decades, asymptotics of the quantization errors have been
studied both for absolutely continuous distributions and for several classes of fractal measures , including self-similar measures, $F$-conformal measures, self-affine measures,
 Markov-type measures and some in-homogeneous self-similar measures (cf. \cite{GL:01,GL:02,GL:08,GL:12,KZ:15,KZ:16,LM:02,PK:04,RM:12,Zhu:08b,Zhu:13,Zhu:18}).

 Recently, Kesseb\"{o}hmer et al proved a surprising general result that the upper quantization dimension for an arbitrary compactly supported probability measures agrees with its R\'{e}nyi dimension at the point $q_r$ where the $L_q$-spectrum and the line through the origin with slope $r$ intersect (see \cite{KNZ:22}). This result, along with some results in \cite{NX:19,PS:00}, leads to the existence of the quantization dimensions for arbitrary self-conformal measures and some nice formulas for the exact values. Further, this work will certainly be helpful for us to establish some explicit and computable formulas for various classes of fractal measures, such as self-affine measures as studied in \cite{Feng:05}.

In the above-mentioned research on the quantization for fractal measures, certain kind of separation properties, such as the open set condition (OSC) or the strong separation condition, are required. As a consequence, the measures for certain sets can be well estimated in a convenient manner. Up to now, little is known about the asymptotics of the quantization errors for those measures with overlaps, where the measures exhibit very complicated structure and can hardly be estimated even for those "good" sets.

In the present paper, we study the asymptotics of the quantization errors for the image measures of Markov-type measures with complete overlaps associated with a class of graph-directed IFS. We emphasize that we are concerned with the \emph{overlaps for measures}, instead of that for the structure of their support, although, in general, the support can be very complicated and may fail to be a Mauldin-Williams fractal in the normal sense. Due to the complete overlaps, we need only a directed graph of $N$ vertices. The overlaps for measures here are induced by the larger directed graph consisting of $2N$ vertices which possesses abundant possibility of admissible paths. Furthermore,  we will construct examples to show that, the image measures of the Markov-type measures are typically no longer of Markov-type, and in general, they are not equivalent to Markov-type measures.

Our work is partially motivated by our previous study on a class of inhomogenous self-similar measures $\nu$ supported on self-similar sets (cf. \cite{Zhu:16} and Example \ref{eg1}). These measures turn out to be extreme examples of the present work. We note that Roychowdhury (\cite{RM:12}) made efforts and obtained some bounds for the upper and lower quantization dimension for the above-mentioned extreme cases. However, he was unable to prove the existence of the quantization dimension, even with an extremely strong assumption as stated in
\cite[Remark 3.7]{RM:12}. In \cite{Zhu:16}, we proved that the quantization dimension for this class of measures exists and established an explicit formula for its exact value which is really computable and established a sufficient condition for the upper quantization coefficient to be finite. An easy corollary of Theorem \ref{mthm} of the present paper will provide a necessary and sufficient condition for this.
\subsection{The quantization error and its asymptotics}
Let $r\in (0,\infty)$ and $k\in\mathbb{N}$. Let $d$ denote the Euclidean metric on $\mathbb{R}^q$. For every $k\geq 1$, let $\mathcal{D}_k:=\{\alpha\subset\mathbb{R}^q:1\leq {\rm
card}(\alpha)\leq k\}$. For $x\in\mathbb{R}^q$ and $\emptyset\neq A,B\subset\mathbb{R}^q$, let $d(x,A):=\inf_{a\in A}d(x,a)$ and $d(A,B):=\inf_{a\in A,b\in B}d(a,b)$. The $k$th quantization
error for $\nu$ of order $r$ can be defined by
\begin{eqnarray}\label{quanerror}
e^r_{k,r}(\nu)=\inf\limits_{\alpha\in\mathcal{D}_k}\int d(x,\alpha)^{r}d\nu(x).
\end{eqnarray}
One may see \cite{GL:00} for various equivalent definitions of the quantization error and interesting interpretations in different contexts. For $r\in[1,\infty)$, $e_{k,r}(\nu)$ is equal to
the minimum error when approximating $\nu$ by discrete probability measures supported on at most $k$ points in the $L_r$-metrics.

The asymptotics of the quantization errors can be characterized by the upper and lower quantization coefficient. For $s\in (0,\infty)$, the
$s$-dimensional upper and lower quantization coefficient for $\nu$ of order $r$ can be defined by
\[
\overline{Q}_r^s(\nu):=\limsup_{k\to\infty}k^{\frac{r}{s}}e_{k,r}^r(\nu),\;\underline{Q}_r^s(\nu):=\liminf_{k\to\infty}k^{\frac{r}{s}}e_{k,r}^r(\nu).
\]
The upper (lower) quantization dimension for $\nu$ of order $r$ is exactly the critical point at which the upper (lower) quantization
coefficient jumps from zero to infinity (cf. \cite{GL:00,PK:01}), which are given by
\[
\overline{D}_r(\nu)=\limsup_{k\to\infty}\frac{\log k}{-\log e_{k,r}(\nu)};\;\;\underline{D}_r(\nu)=\liminf_{k\to\infty}\frac{\log k}{-\log e_{k,r}(\nu)}.
\]
If $\overline{D}_r(\nu)=\underline{D}_r(\nu)$, then we denote the common value by $D_r(\nu)$. Compared with the estimate for the upper and lower quantization dimension, it is certainly much more difficult to examine the finiteness and positivity of the upper and lower quantization coefficient. This is similar to the relationship between the Hausdorff dimension of a set and its Hausdorff measure in the exact dimension.

Next, let us recall a classical result by Graf and Luschgy. Let $(f_i)_{i=1}^N$ be a family of contractive similarity mappings on $\mathbb{R}^q$ with contraction ratios $(c_i)_{i=1}^N$. According to \cite{Hut:81}, there exists a unique non-empty compact set $E$ satisfying the equation
$E=\bigcup_{i=1}^Nf_i(E)$.
 The set $E$ is called the self-similar set determined by  $(f_i)_{i=1}^N$. Given a probability vector $(p_i)_{i=1}^N$, there exists a unique Borel probability measure satisfying
$ \nu=\sum_{i=1}^Np_i\nu\circ f_i^{-1}$. We call $\nu$ the self-similar measure associated with $(f_i)_{i=1}^N$ and $(p_i)_{i=1}^N$.
 We say that $(f_i)_{i=1}^N$ satisfies the OSC if there exists a bounded non-empty open set $U$ such that $f_i(U),1\leq i\leq N$, are pairwise disjoint and $f_i(U)\subset U$ for all
 $1\leq i\leq N$. Let $\xi_r$ be implicitly defined by $\sum_{i=1}^N(p_ic_i^r)^{\frac{\xi_r}{\xi_r+r}}=0$.
 Assuming the OSC for $(f_i)_{i=1}^N$, Graf and Luschgy proved that $0<\underline{Q}_r^{\xi_r}(\nu)\leq\overline{Q}_r^{\xi_r}(\nu)<\infty$ (cf. \cite{GL:01,GL:02}). This work has enlightened almost all subsequent study on the quantization for fractal measures.

\subsection{A class of Mauldin-Williams fractals}
Mauldin-Williams (MW) fractals were introduced and studied in detail in \cite{MW:88}. Multi-fractal decompositions for such fractals were accomplished by Edgar and Mauldin (cf.
\cite{EM:92}). Next, we describe a class of MW-fractals with complete overlaps, which we will work with in the remainder of the paper.

 Fix an integer $N\geq 2$. Let $\mathcal{G}$ be a directed graph  with vertices $1,2,\ldots, 2N$. We assume that there exists at most one directed edge from a vertex $i$ to another vertex
 $j$, and that there exists at least one edge leaving a vertex $i$. Let $(c_{i,j})_{i,j=1}^{2N}$ be an incidence matrix for $\mathcal{G}$. Thus, $c_{i,j}\geq 0$ for all $1\leq i,j\leq 2N$, and $c_{i,j}>0$ if and only if there exists one directed
 edge from the vertex $i$ to $j$. We write
\begin{eqnarray*}
&&G_1=\Omega:=\{1,\ldots, 2N\};\;\;\Psi:=\{1,\ldots, N\};\;\;\Psi^*:=\bigcup_{n=1}^\infty\Psi^n;\;\;\Omega^*:=\bigcup_{n=1}^\infty\Omega^n.\\
&&G_k:=\{(\sigma_1,\ldots,\sigma_k)\in\Omega^k:c_{\sigma_i\sigma_{i+1}}>0\;\;1\leq i\leq k-1\};\;\;k\geq
2;\\&&G_\infty:=\{(\sigma_1,\ldots,\sigma_k,\ldots)\in\Omega^{\mathbb{N}}:c_{\sigma_i\sigma_{i+1}}>0,\;\;i\geq 1\};\;G^*:=\bigcup_{k=1}^\infty G_k.
\end{eqnarray*}
For $i\in\Psi$, we define $i^+:=i+N$ and for $(i,j)\in\Psi^2$, we define
\begin{equation}\label{nj}
\mathcal{N}_{i,j}:=\{(i,j),(i,j^+),(i^+,j),(i^+,j^+)\};\;\;\mathcal{M}_{i,j}:=\mathcal{N}_{i,j}\cap G_2.
\end{equation}

Let $J_i, i\in\Psi$, be non-empty, pairwise disjoint compact subsets of $\mathbb{R}^q$ with $\overline{{\rm int}(J_i)}=J_i, i\in\Psi$, where $\overline{A}$ and ${\rm int}(A)$ denote the
closure and interior of a subset $A$ of $\mathbb{R}^q$. We define $J_{i^+}:=J_i$ for $i\in\Psi$. For every $\sigma\in\Psi^*\cup\Psi^{\mathbb{N}}$, we define
\begin{equation}\label{zs1}
\mathcal{T}(\sigma):=\left\{\begin{array}{ll}
\big\{\widetilde{\sigma}\in\Omega^*:\;\widetilde{\sigma}_i=\sigma_i,\;{\rm or}\;\sigma_i^+, 1\leq i\leq n\big\}&{\rm if}\;\sigma\in\Psi^n\\
\big\{\widetilde{\sigma}\in\Omega^{\mathbb{N}}:\;\widetilde{\sigma}_i=\sigma_i,\;{\rm or}\;\sigma_i^+,i\geq 1\big\}&{\rm if}\;\sigma\in\Psi^{\mathbb{N}}
\end{array}\right..
\end{equation}
We will use the following sets in the construction of MW-fractals with overlaps:
\begin{eqnarray}\label{sk}
&&\mathcal{S}_1:=\Psi,\;\mathcal{S}_k=\{\sigma\in\Psi^k: \mathcal{T}(\sigma)\cap G_k\neq\emptyset\},\;k\geq 2;\\
&&\mathcal{S}_\infty=\{\sigma\in\Psi^{\mathbb{N}}: \mathcal{T}(\sigma)\cap G_\infty\neq\emptyset\};\;\mathcal{S}^*=\bigcup_{k\geq 1}\mathcal{S}_k.\nonumber
\end{eqnarray}

Let $T_{i,j},(i,j)\in \mathcal{S}_2$, be contractive similarity mappings on $\mathbb{R}^q$ of similarity ratios $s_{i,j},(i,j)\in \mathcal{S}_2$. Note that for every
$(\widetilde{i},\widetilde{j})\in G_2$, there exists a unique $(i,j)\in\mathcal{S}_2$ such that $(\widetilde{i},\widetilde{j})\in \mathcal{M}_{i,j}$. We assume that
\begin{equation}\label{overlapping mapps}
T_{\widetilde{i},\widetilde{j}}:=T_{i,j},\;\;{\rm for\;every}\;\; (\widetilde{i},\widetilde{j})\in \mathcal{M}_{i,j}.
\end{equation}
Thus, we obtain an IFS $\mathcal{T}=\{T_{i,j}: (i,j)\in G_2\}$.  $\mathcal{T}$ has \emph{complete overlaps} if for some $(i,j)\in \Psi^2$, ${\rm card}(M_{i,j})\geq 2$ ; for example, if
$(i,j), (i,j^+)\in G_2$, then
\[
J_{(i,j^+)}:=T_{i,j^+}(J_{j^+})=T_{i,j}(J_j)=J_{(i,j)}.
\]
We just write $J_{i,j}$ for $J_{(i,j)}$. Now we further assume that, for every $i\in\Psi$,
\begin{equation}\label{ssc}
\bigcup_{j:(i,j)\in\mathcal{S}_2}T_{i,j}(J_j)\subset J_i;\;\; T_{i,j}(J_j), (i,j)\in\mathcal{S}_2, \;{\rm are\; pairwise\; disjoint}.
\end{equation}
For each $\sigma=(\sigma_1,\ldots, \sigma_n)\in\mathcal{S}_n\cup G_n$, we define
\begin{eqnarray}\label{gamma-sigma}
T_\sigma=\left\{\begin{array}{ll}
id_{\mathbb{R}^q},&{\rm if}\;n=1\\
T_{\sigma_1,\sigma_2}\circ T_{\sigma_2,\sigma_3}\circ\cdots T_{\sigma_{n-1},\sigma_n},&{\rm if}\;n>1
\end{array}\right.;\;\;J_\sigma=T_\sigma(J_{\sigma_n}).
\end{eqnarray}
For every $n\geq 1$, we call the sets $J_\sigma,\sigma\in \mathcal{S}_n\cup G_n$, \emph{cylinders of order $n$}.
Then we obtain a MW-fractal $K$ (\cite[Theorem 1]{MW:88}):
\begin{equation*}\label{MWfractal}
K=\bigcap_{k\geq 1}\bigcup_{\sigma\in G_k}J_\sigma=\bigcap_{k\geq 1}\bigcup_{\sigma\in \mathcal{S}_k}J_\sigma.
\end{equation*}
For and $\widetilde{\sigma}\in G_n\cup\mathcal{S}_n$, we define
\begin{eqnarray*}
s_{\widetilde{\sigma}}:=\left\{\begin{array}{ll}
1,&{\rm if}\;n=1\\
s_{\widetilde{\sigma}_1\widetilde{\sigma}_2}\cdots s_{\widetilde{\sigma}_{k-1}\widetilde{\sigma}_k}&{\rm if}\;n>1
\end{array}\right..
\end{eqnarray*}
\begin{remark}{\rm
In general, the sets $\mathcal{S}_n, n\geq 1$, can be very complicated and the fractal set $K$ may not be generated by a reduced graph-directed IFS (see Example \ref{eg2}). We will impose some conditions for the incidence matrix such that $\mathcal{S}_n$ can be well tracked and $K$ remains a typical MW-fractal. We will then focus on the overlaps for measures induced by the larger transition matrix $P$.
}\end{remark}
\subsection{Markov-type measures with complete overlaps}
Let $\theta$ denote the empty word. For every $\widetilde{\sigma}\in G^*\cup\mathcal{S}^*$, we denote by $|\widetilde{\sigma}|$ the length of $\widetilde{\sigma}$; we define
$|\widetilde{\sigma}|:=\infty$ if $\widetilde{\sigma}\in G_\infty\cup\mathcal{S}_\infty$. For $\widetilde{\sigma}\in G^*\cup G_\infty$, or $\widetilde{\sigma}\in\mathcal{S}^*\cup
\mathcal{S}_\infty$ and $1\leq h\leq |\widetilde{\sigma}|$, we write
\[
\widetilde{\sigma}|_h:=(\widetilde{\sigma}_1,\ldots,\widetilde{\sigma}_h);\;\;{\widetilde{\sigma}^\flat}:=\left\{\begin{array}{ll}
\theta,&{\rm if}\;|\widetilde{\sigma}|=1\\
\widetilde{\sigma}|_{|\widetilde{\sigma}|-1},&{\rm if}\;|\widetilde{\sigma}|>1
\end{array}\right..
\]
For every $h\geq 1$, we call $\widetilde{\sigma}$ a \emph{descendant} of $\widetilde{\sigma}|_h$, and $\widetilde{\sigma}|_h$ a \emph{predecessor} of $\widetilde{\sigma}$.

Let $\pi:G_\infty\to K$ be defined by $\pi(\widetilde{\sigma}):=\bigcap_{k=1}^\infty J_{\widetilde{\sigma}|_k}$.
We define
\[
\Gamma(\sigma):=G_n\cap\mathcal{T}(\sigma),\;\sigma\in\mathcal{S}_n;\;\;[\widetilde{\sigma}]:=\{\widetilde{\tau}\in G_\infty: \widetilde{\tau}|_n=\widetilde{\sigma}\},\;\widetilde{\sigma}\in G^*.
\]
Then for every $\widetilde{\sigma}\in\Gamma(\sigma)$, by (\ref{overlapping mapps}), we have $J_{\widetilde{\sigma}}=J_\sigma$.
By (\ref{ssc}), we have
\begin{eqnarray*}
\pi^{-1}(J_\sigma)=\pi^{-1}(J_\sigma\cap K)=\bigcup_{\widetilde{\sigma}\in\Gamma(\sigma)}[\widetilde{\sigma}],\;\sigma\in\mathcal{S}^*.
\end{eqnarray*}

Let $(\chi_i)_{i=1}^{2N}$ be a positive probability vector. Let $P=(p_{i,j})_{i,j=1}^{2N}$ be a row-stochastic matrix (transition matrix), with $p_{i,j}>0$ if and only if $(i,j)\in G_2$.
For every $n\geq 1$ and $\widetilde{\sigma}\in G_n$, we define
\begin{eqnarray*}
p_{\widetilde{\sigma}}:=\left\{\begin{array}{ll}
1,&{\rm if}\;n=1\\
p_{\widetilde{\sigma}_1\widetilde{\sigma}_2}\cdots p_{\widetilde{\sigma}_{n-1}\widetilde{\sigma}_n}&{\rm if}\;n>1
\end{array}\right..
\end{eqnarray*}
By Kolmogorov consistency theorem, there exists a unique Borel probability measure $\nu$ on $G_\infty$ such that for every $k\geq 1$ and $\widetilde{\sigma}=(\widetilde{\sigma}_1,\ldots,
\widetilde{\sigma}_k)\in G_k$,
\begin{equation}\label{nu}
\nu([\widetilde{\sigma}])=\chi_{\widetilde{\sigma}_1}p_{\widetilde{\sigma}}=\chi_{\widetilde{\sigma}_1}p_{\widetilde{\sigma}_1\widetilde{\sigma}_2}\cdots
p_{\widetilde{\sigma}_{k-1}\widetilde{\sigma}_k}.
\end{equation}
We define $\mu=\nu\circ \pi^{-1}$. We call $\mu$ a \emph{Markov-type measure with complete overlaps} if ${\rm card}(M_{i,j})\geq 2$ for some $(i,j)\in\mathcal{S}_2$. We have
\begin{equation}\label{mu-meausre}
\mu(J_\sigma)=\nu\circ \pi^{-1}(J_\sigma)=\sum_{\widetilde{\sigma}\in\Gamma(\sigma)}\chi_{\widetilde{\sigma}_1}p_{\widetilde{\sigma}},\;\sigma\in\mathcal{S}^*.
\end{equation}

\subsection{Statement of the main results}
Let $P$ be the $2N\times 2N$ transition matrix as above. We write
\begin{eqnarray*}
P=\left(\begin{array}{cccc}
P_1 & P_3\\
P_4 &   P_2 \\
\end{array}\right)\; {\rm with }\;\;P_1=\big(p_{i,j}\big)_{i,j=1}^N.
\end{eqnarray*}
Let $\textbf{0}$ denote a zero matrix. In the present paper, we consider two cases:

\textbf{Case I:} $P$ is reducible. We assume:

(A1) $P_1, P_2$ are both irreducible and $P_4=\textbf{0}$;

(A2) ${\rm card}(\{j\in\Psi:p_{i,j}>0\}), {\rm card}(\{j\in\Psi:p_{i^+,j^+}>0\})\geq 2$; $i\in\Psi$;

(A3) for $(i,j)\in\Psi^2$, either $\mathcal{M}_{ij}=\emptyset$, or
$\mathcal{M}_{ij}=\{(i,j),(i,j^+),(i^+,j^+)\}$.

\textbf{Case II:} $P$ is irreducible. We will assume (A2) and

(A4) for $(i,j)\in\Psi^2$, either $\mathcal{M}_{i,j}=\emptyset$ or $\mathcal{M}_{i,j}=\mathcal{N}_{i,j}$ (see (\ref{nj})).

(A5) $P_1$ is irreducible.

For every $s\in (0,\infty)$, we define
\begin{eqnarray*}
&A_1(s):=\big((p_{i,j}s_{i,j}^r)^s\big)_{i,j=1}^N;\;\;\;A_2(s):=\big((p_{i^+,j^+}s_{i^+,j^+}^r)^s\big)_{i,j=1}^N;\\
&A_3(s):=\big((p_{i,j^+}s_{i,j^+}^r)^s\big)_{i,j=1}^N;\;\;\;A_4(s):=\big((p_{i^+,j}s_{i^+,j}^r)^s\big)_{i,j=1}^N.
\end{eqnarray*}
Let $\psi_i(s)$ denote the spectral radius of $A_i(s)$ and $\rho_i(s):=\psi_i\big(\frac{s}{s+r}\big)$. Define
\begin{eqnarray*}\label{A(s)}
A(s)=\left(\begin{array}{cc}
A_1(s) & A_3(s)\\
A_4(s) &   A_2(s)) \\
\end{array}\right).
\end{eqnarray*}
Let $\psi(s)$ denote the spectral radius of $A(s)$ and define $\rho(s):=\psi\big(\frac{s}{s+r}\big)$.

\begin{remark}
{\rm By \cite[Theorem 2]{MW:88}, $\psi_i(s),\psi(s)$ are continuous and strictly decreasing.
Assuming (A2), we have $\psi_i(0)\geq 2$ and $\psi_i(1)<1$ for $i=1,2$. Thus, there exists a unique $s\in(0,1)$ such that $\psi_i(s)=1$. Further,
there exists a unique positive number
$s_{i,r}$ such that $\rho_i(s_{i,r})=1$. Also, there exists a unique $s_r$ with $\rho(s_r)=1$. When $P_4=\textbf{0}$, we have $s_r=\max\{s_{1,r},s_{2,r}\}$.
}
\end{remark}

Let $\pi_1:\mathcal{S}_\infty\mapsto K$ be defined by
$\pi_1(\sigma):=\bigcap_{k=1}^\infty J_{\sigma|_k}$.
By (\ref{ssc}), $\pi_1$ is a bijection. We say that $\mu$ \emph{is reducible} if $\mu=\nu_1\circ\pi_1^{-1}$ for the Markov-type measure $\nu_1$ associated with
some transition matrix $\widetilde{P}_{N\times N}$ and some initial probability vector $(\widetilde{\chi}_i)_{i=1}^N$.  Whenever $\mu$ is reducible, the asymptotics
of $(e_{n,r}(\mu))_{n=1}^\infty$ is characterized by \cite[Theorem 1.1]{KZ:15}. For every $i\in\Psi$, we define
\begin{eqnarray*}
\mathcal{S}_n(i):=\{\sigma\in \mathcal{S}_n:\sigma_n=i\},\;n\geq 1;\;\;\mathcal{S}^*(i):=\bigcup_{n=1}^\infty\mathcal{S}_n(i).
\end{eqnarray*}
For every $n\geq 1$ and $\sigma\in\mathcal{S}_n(i)$, we split $\mu(J_\sigma)$ into two parts (see (\ref{mu-meausre})):
\begin{eqnarray*}
I_{1,\sigma}:=\sum_{\widetilde{\sigma}\in\Gamma(\sigma),\widetilde{\sigma}_n=i}\chi_{\widetilde{\sigma}_1}p_{\widetilde{\sigma}};\;\;
I_{2,\sigma}:=\sum_{\widetilde{\sigma}\in\Gamma(\sigma),\widetilde{\sigma}_n=i^+}\chi_{\widetilde{\sigma}_1}p_{\widetilde{\sigma}}.
\end{eqnarray*}

In the remaining part of the paper, we always assume that the IFS $\mathcal{T}$ satisfies (\ref{overlapping mapps}) and (\ref{ssc}), and $\mu$ always denotes the measure as defined in (\ref{mu-meausre}). As our first result, we provide necessary and sufficient conditions for $\mu$ to be reducible, with the assumption of (A2) and (A5), or (A1)-(A3).

\begin{theorem}\label{mthm1}
(1) Assume that (A2) and (A5) hold. Then $\mu$ is reducible if and only if for every $i\in\Psi$, either one of the following holds:
\begin{itemize}
\item[(a)] $p_{i,j}+p_{i,j^+}=p_{i^+,j}+p_{i^+,j^+}$ for every $j\in\Psi$;
\item[(b)] $\chi_{i^+}I_{1,\sigma}=\chi_iI_{2,\sigma}$ for every $\sigma\in\mathcal{S}^*(i)$.
\end{itemize}
(2) Assume that (A1)-(A3) hold. Then $\mu$ is reducible if and only if for every $(i,j)\in \mathcal{S}_2$, we have
$p_{i,j}+p_{i,j^+}=p_{i^+,j^+}$.
\end{theorem}
As our second result, for Case I, we prove that $D_r(\mu)=s_r$ and
establish a necessary and sufficient condition for $\underline{Q}_r^{s_r}(\mu)$ and $\overline{Q}_r^{s_r}(\mu)$ to be both positive and finite. That is,
\begin{theorem}\label{mthm}
Assume that (A1)-(A3) hold. We have
\begin{enumerate}
\item[\rm (i)]$D_r(\mu)=s_r,\;\underline{Q}_r^{s_r}(\mu)>0$;
\item[\rm (ii)] $\overline{Q}_r^{s_r}(\mu)<\infty$ if and only if $s_{1,r}\neq s_{2,r}$;
\item[\rm (iii)]there exists some $r_0>0$, such that $\overline{Q}_r^{s_r}(\mu)<\infty$ for every $r\in(0,r_0)$.
\end{enumerate}
\end{theorem}

For the proof of Theorem \ref{mthm}, we will construct some
 auxiliary measures by applying some ideas of Mauldin and Williams \cite{MW:88}. These
measures will allow us to estimate the quantization error for $\mu$ in a more accurate and concise way.  It seems somewhat surprising that Theorem \ref{mthm} (i) and (ii) are shared by Markov-type measures in non-overlapping cases (see \cite{KZ:15}). When $P$ is irreducible, $\mathcal{T}$ is even more overlapped. Our third result shows that the measure $\mu$ exhibits quite different properties and Theorem \ref{mthm} can fail. Assuming (A2), (A4) and (A5), we will show
 that there exists a unique positive number $t_r$ satisfying
\begin{equation}\label{zsg4}
\lim_{n\to\infty}\frac{1}{n}\log\sum_{\sigma\in\mathcal{S}_n}\bigg(\sum_{\widetilde{\sigma}\in\Gamma(\sigma)}
\chi_{\widetilde{\sigma}_1}p_{\widetilde{\sigma}}s_{\widetilde{\sigma}}^r\bigg)^{\frac{t_r}{t_r+r}}=0.
\end{equation}
We will consider the following two cases which might help to illustrate Case II:
\begin{enumerate}
\item[\rm(g1)] $p_{i,j}+p_{i,j^+}=p_{i^+,j}+p_{i^+,j^+}$ for every $(i,j)\in\mathcal{S}_2$;

\item[\rm (g2)] $p_{i,j}+p_{i^+,j}=p_{i,j^+}+p_{i^+,j^+}$ for every $(i,j)\in\mathcal{S}_2$.
\end{enumerate}
\begin{theorem}\label{mthm2}
Assume that (A2), (A4) and (A5) hold. Then we have
\begin{enumerate}
\item[\rm (i)]$D_r(\mu)=t_r\leq s_r$, and $0<\underline{Q}_r^{t_r}(\mu)\leq\overline{Q}_r^{t_r}(\mu)<\infty$;
\item[\rm (ii)] if (g1) or (g2) holds, then we have $t_r<s_r$.
\end{enumerate}
\end{theorem}

For the proof of Theorem \ref{mthm2}, we will apply a Helley-type theorem (cf. \cite[Theorem 1.23]{PM:95}) and
some ideas contained in the proof of \cite[Theorem 5.1]{Fal:97} to construct some auxiliary measures. These measures are closely connected with the upper and lower quantization coefficient for $\mu$ and enable us to prove Theorem \ref{mthm2} (i) in a convenient way.

\section{Proof of Theorem \ref{mthm1} and some examples}

Let $\Psi=\{1,2,\ldots, N\}$ as before. For each $n\geq 1$ and $i\in\Psi$, we write
\begin{eqnarray*}
&&H_{1,n}:=G_n\cap\Psi^n;\;H_1^*:=G^*\cap\Psi^*;\;\;H_1^\infty:=G_\infty\cap\Psi^{\mathbb{N}};\\
&&H_{1,n}(i):=\{\sigma\in H_{1,n}:\sigma_1=i\};\;\;H_1^*(i):=\{\sigma\in H_1^*:\sigma_1=i\};\\&&H_1^\infty(i):=\{\sigma\in H_1^\infty:\sigma_1=i\}.
\end{eqnarray*}
Let $\Upsilon:=\{1^+,2^+,\ldots, N^+\}$. Let $H_{2,n}, H_2^*, H_2^\infty, H_{2,n}(i^+), H_2^*(i^+), H_2^\infty(i^+)$, be defined in the same manner by replacing $\Psi$ with $\Upsilon$.

For $\sigma,\tau\in G^*$ or $\sigma,\tau\in\mathcal{S}^*$, we denote by $\sigma\ast\tau$ the concatenation of $\sigma$ and $\tau$.  For every $\sigma\in\Psi^*$, we write
$\sigma^+:=(\sigma_1^+,\ldots,\sigma_{|\sigma|}^+)$. We define
\begin{eqnarray*}
&&\mathcal{L}(\sigma):=\{\sigma, \sigma^+\}\cup \{\sigma|_h\ast(\sigma_{h+1}^+,\sigma_{h+2}^+,\ldots,\sigma_n^+),1\leq h\leq n-1\};\;\sigma\in\mathcal{S}_n.
\end{eqnarray*}

The subsequent two lemmas shows that with (A1)-(A3), or (A2) and (A4), the sets $\Gamma(\sigma)$ and $\mathcal{S}_n$ can be well tracked.
\begin{lemma}\label{condition1}
 Assume that (A1)-(A3) hold. We have
 \begin{eqnarray*}\label{tem3}
 \Gamma(\sigma)=\mathcal{L}(\sigma),\;\sigma\in \mathcal{S}_n;\;\mathcal{S}_n=H_{1,n},\;n\geq 1;\;\mathcal{S}^*=H_1^*.
 \end{eqnarray*}
 \end{lemma}
 \begin{proof}
 Let $\sigma\in\mathcal{S}_n$ and $\widetilde{\sigma}\in\Gamma(\sigma)$ be given. By (A1), for every $1\leq i\leq n-1$, we have $p_{\sigma_i^+,\sigma_{i+1}}=0$. Hence, if
 $\widetilde{\sigma}_{h+1}=\sigma_{h+1}^+$ for some integer $h\geq 0$, then $\widetilde{\sigma}_l=\sigma_l^+$ for all $h+1\leq l\leq n$. This implies that either $\widetilde{\sigma}=\sigma$
 or $\sigma^+$, or for some $1\leq h\leq n-1$,
$\widetilde{\sigma}=\sigma|_h\ast(\sigma_{h+1}^+,\sigma_{h+2}^+,\ldots,\sigma_n^+)$. It follows that $\Gamma(\sigma)\subset\mathcal{L}(\sigma)$.

 For every $\sigma\in\mathcal{S}_n$, by (\ref{sk}) and (A1), we have $\sigma\in G_n$, or $\sigma^+\in G_n$, or for some $1\leq h\leq n-1$,
we have $\sigma|_h\ast(\sigma_{h+1}^+,\sigma_{h+2}^+,\ldots,\sigma_n^+)\in G_n$.
 Using this and (A3), we obtain that $ \mathcal{L}(\sigma)\subset G_n$.
It follows that $ \mathcal{L}(\sigma)\subset \Gamma(\sigma)$.
Thus the first part of the lemma holds. In particular, we have $\sigma\in G_n$, which implies that $\mathcal{S}_n\subset H_{1,n}$.
Since $H_{1,n}\subset\mathcal{S}_n$, we obtain that $\mathcal{S}_n=H_{1,n}$ and $\mathcal{S}^*=H_1^*$.
 \end{proof}

Let $\mathcal{T}(\sigma)$ be as defined in (\ref{zs1}). We have
\begin{lemma}\label{T5}
 Assume that (A2) and (A4) hold. Then
\begin{eqnarray*}
\Gamma(\sigma)=\mathcal{T}(\sigma),\;\sigma\in \mathcal{S}_n;\;\mathcal{S}_n=H_{1,n},\;n\geq 1;\;\mathcal{S}^*=H_1^*.
\end{eqnarray*}
\end{lemma}
\begin{proof}
This can be proved analogously to the proof of Lemma \ref{condition1}.
\end{proof}

 Let $n\geq 1$ and $i\in\Psi$, we write
\[
 B_n(i):=\{\sigma\in  H_{1,n}:\sigma_n=i\},\;\; B_n(i^+):=\{\sigma^+\in H_{2,n}:\sigma_n^+=i^+\}.
\]
We clearly have $ B_n(i)\subset\mathcal{S}_n(i)$. If (A1)-(A3) hold, or (A2) and (A4) hold, then by Lemmas \ref{condition1} and \ref{T5}, we have $ B_n(i)=\mathcal{S}_n(i)$.

\begin{remark}\label{Zhu6}{\rm
Assume that (A5) holds. By induction, one can easily see that $ B_n(i)\neq\emptyset$ for every $i\in\Psi$ and $n\geq 1$.
}\end{remark}
\textbf{Proof of Theorem \ref{mthm1} (1)}

For every $(i,j)\in\mathcal{S}_2$, we define
\begin{equation}\label{zsg5}
\widetilde{p}_{i,j}:=\frac{\mu(J_{i,j})}{\mu(J_i)}=\frac{\chi_i(p_{i,j}+p_{i,j^+})+\chi_{i^+}(p_{i^+,j}+p_{i^+,j^+})}{\chi_i+\chi_{i^{+}}}.
\end{equation}
Then we have $\mu(J_{i,j})=\mu(J_i)\cdot \widetilde{p}_{i,j}$ for every $(i,j)\in\mathcal{S}_2$.

 $\Rightarrow)$ First we assume that $\mu=\nu\circ \pi_1^{-1}$ for some markov-type measure $\nu$. Then $\nu$ is the Markov-type measure associated with
 $\widetilde{P}=(\widetilde{p}_{i,j}\}_{i,j=1}^N$ and $\widetilde{\chi}=(\chi_i+\chi_{i^+})_{i=1}^N$. For every $i\in\Psi$ and $\sigma\in\mathcal{S}^*(i)$ and
 $j$ with $(i, j)\in\mathcal{S}_2$, we have
\begin{eqnarray}\label{temp1}
\Delta_{\sigma,j}:=\mu(J_{\sigma\ast j})-\mu(J_\sigma)\cdot\widetilde{p}_{i,j}=0.
\end{eqnarray}
Note that $I_{1,\sigma}+I_{2,\sigma}=\mu(J_\sigma)$. It follows that
\begin{eqnarray}\label{temp2}
\Delta_{\sigma,j}&=&I_{1,\sigma}(p_{i,j}+p_{i,j^+})+I_{2,\sigma}(p_{i^+,j}+p_{i^+,j^+})-\mu(J_\sigma)\cdot\widetilde{p}_{i,j}\nonumber\\
&=&I_{1,\sigma}(p_{i,j}+p_{i,j^+}-\widetilde{p}_{i,j})+I_{2,\sigma}(p_{i^+,j}+p_{i^+,j^+}-\widetilde{p}_{i,j})\nonumber\\
&=&\frac{\chi_{i^+}I_{1,\sigma}-\chi_i I_{2,\sigma}}{\chi_i+\chi_{i^+}}(p_{i,j}+p_{i,j^+}-p_{i^+,j}-p_{i^+,j^+}).
\end{eqnarray}
By (\ref{temp1}) and (\ref{temp2}), for every $i\in\Psi$, we have the following two cases:

Case (1): for every $j$ with $(i,j)\in\mathcal{S}_2$, we have
\begin{equation}\label{sg1}
p_{i,j}+p_{i,j^+}=p_{i^+,j}+p_{i^+,j^+}.
\end{equation}

Case (2): there exists some $j_0\in\Psi$ with $(i,j_0)\in\mathcal{S}_2$ such that (\ref{sg1}) fails. Then by (\ref{temp1}) and (\ref{temp2}), we obtain that
\begin{equation}\label{sg2}
\chi_{i^+}I_{1,\sigma}=\chi_i I_{2,\sigma},\;\;{\rm for\;\;every}\;\;\sigma\in\mathcal{S}^*(i),
\end{equation}
Otherwise, there would be some $\sigma\in\mathcal{S}^*(i)$ such that (\ref{temp1}) fails for $j_0$. This contradicts the assumption that $\mu$ is reducible.

$\Leftarrow)$ Assume that, for every $i\in\Psi$, (a) or (b) holds.
Then from (\ref{temp2}), we know that (\ref{temp1}) holds for every $\sigma\in\mathcal{S}^*(i)$ and every $j$ with $(i,j)\in\mathcal{S}_2$, which implies that $\mu=\nu\circ\pi^{-1}$ for
the Markov-type measure associated with $P=(\widetilde{p}_{i,j})_{i,j=1}^N$ and $\widetilde{\chi}=(\chi_i+\chi_{i^+})_{i=1}^N$.

Next, we give sufficient conditions such that $\mu$ is reducible or non-reducible.
\begin{corollary}\label{cor1}
Assume that (A2), (A4) and (A5) hold and that
\begin{enumerate}
\item[(b1)] $\chi_i=\chi_{i^+}$ for every $i\in\Psi$;
\item[(b2)] for every $i\in\Psi$, there exists some $l\in\Psi$ such that $p_{l,i}+p_{l^+,i}\neq p_{l,i^+}+p_{l^+,i^+}$.
\end{enumerate}
Then $\mu$ is reducible if and only if (\ref{sg1}) holds for every $(i,j)\in\mathcal{S}_2$.
\end{corollary}
\begin{proof}
By (b1) and (b2), for $\sigma:=(l,i)\in\mathcal{S}_2(i)$, we have
\begin{equation}\label{z7}
\frac{I_{1,\sigma}}{I_{2,\sigma}}=\frac{\chi_lp_{l,i}+\chi_{l^+}p_{l^+,i}}{\chi_lp_{l,i^+}+\chi_{l^+}p_{l^+,i^+}}
=\frac{p_{l,i}+p_{l^+,i}}{p_{l,i^+}+p_{l^+,i^+}}\neq 1=\frac{\chi_i}{\chi_{i^+}}.
\end{equation}
The corollary follows immediately from Theorem \ref{mthm1} (1).
\end{proof}
\begin{corollary}
Assume that (A2), (A4), (A5) and (b1) hold and that for every $(l,i)\in\mathcal{S}_2$, we have
 $p_{l,i}+p_{l^+,i}=p_{l,i^+}+p_{l^+,i^+}$.  Then $\mu$ is reducible.
\end{corollary}
\begin{proof}
We first will show that

\emph{Claim 1}: $\chi_{i^+}I_{1,\sigma}=\chi_i I_{2,\sigma}$ holds for every $i\in\Psi$, every $n\geq1$ and $\sigma\in\mathcal{S}_n(i)$.
 By the hypothesis, we know that $\chi_l=\chi_{l^+}$ for all $l\in\Psi$. Also, for $n=1$, we have $\mathcal{S}_1(i)=\{i\}$. For $\sigma=i\in\mathcal{S}_1(i)$, we have
 $I_{1,\sigma}=\chi_i=\chi_{i^+}=I_{2,\sigma}$. For $n=2$ and every $\sigma=(l,i)\in\mathcal{S}_2(i)$, by considering (\ref{z7}), we have $I_{1,\sigma}=I_{2,\sigma}$. Thus, Claim 1 holds
 for every $i\in\Psi$ and $n=1,2$, and every $\sigma\in\mathcal{S}_n(i)$.

Now we assume that Claim 1 holds for $n=k\geq 2$ and every $i\in\Psi$ and $\sigma\in\mathcal{S}_k(i)$. Let $n=k+1, i\in\Psi$ and $\sigma\in\mathcal{S}_n(i)$. Let $\tau:=\sigma^\flat$. By
Lemma \ref{T5},
\[
\Gamma(\sigma)=\{\widetilde{\tau}\ast i,\; \widetilde{\tau}\ast i^+:\widetilde{\tau}\in\Gamma(\tau)\}.
\]
Note that $\tau\in \mathcal{S}_k(\tau_k)$. By the inductive assumption, we have $I_{1,\tau}=I_{2,\tau}$.
Using the hypothesis of the corollary, we deduce
\begin{eqnarray*}
I_{1,\sigma}&=&\sum_{\widetilde{\tau}\in\Gamma(\tau)}\chi_{\widetilde{\tau}_1}p_{\widetilde{\tau}}\cdot p_{\widetilde{\tau}_k,i}
=I_{1,\tau}p_{\tau_k,i}+I_{2,\tau}p_{\tau_k^+,i}\\
&=&I_{1,\tau}\big(p_{\tau_k,i}+p_{\tau_k^+,i}\big)
=I_{1,\tau}\big(p_{\tau_k,i^+}+p_{\tau_k^+,i^+}\big)=I_{2,\sigma}.
\end{eqnarray*}
By induction, Claim 1 holds. Thus, by Theorem \ref{mthm} (1),  $\mu$ is reducible.
\end{proof}

For two variables $X,Y$ taking values in $(0,\infty)$, we write $X\lesssim Y$ ($X\gtrsim Y$), if there exists some constant $C>0$, such that the inequality $X\leq C Y$ ($X\geq CY$) always holds.
We write $X\asymp Y$, if we have both $X\lesssim Y$ and $X\gtrsim Y$. We define
\[
\underline{\chi}:=\min\limits_{1\leq i\leq 2N}\chi_i,\;\;\overline{\chi}:=\max\limits_{1\leq i\leq 2N}\chi_i.
\]
To complete the proof for Theorem \ref{mthm1} (2), we need one more lemma.

\begin{lemma}\label{condition3}
Assume that $(A1)-(A3)$ hold. For every $i\in\Psi$, there exists some integer $k_0$ such that for every $n\geq k_0$ and some $\sigma\in\mathcal{S}_n(i)$ such that
\begin{equation}\label{s4}
I_{1,\sigma}<\chi_i\chi_{i^+}^{-1} I_{2,\sigma}.
\end{equation}
\end{lemma}
\begin{proof}
Let $\rho_1,\rho_2$ denote the spectral radius of $P_1, P_2$. We first show that, for every $i\in\Psi$, we have
\begin{equation}\label{s3}
\sum_{\sigma\in B_n(i)}p_\sigma\asymp\rho_1^{n-1},\;\;\sum_{\sigma^+\in B_n(i^+)}p_{\sigma^+}\asymp 1.
\end{equation}
By the assumption (A3), we see that $\sum_{j=1}^Np_{i,j}<1$ for every $i\in\Psi$.
By \cite[Theorem 8.1.22]{HJ:13}, this implies that $\rho_1<1$. By (A1), we have $P_4=\textbf{0}$. Thus, $\sum_{j=1}^Np_{i^+,j^+}=1$, for every $i\in\Psi$, which implies that $\rho_2=1$. Let
$(c_{j,i})_{j=1}^N$  and $(c^+_{j,i})_{j=1}^N$ denote the $i$th column of $P_1^{n-1}$ and $P_2^{n-1}$ respectively. Then
\begin{eqnarray*}\label{tem4}
\sum_{\sigma\in B_n(i)}p_\sigma=\sum_{j=1}^Nc_{j,i};\;\;\sum_{\sigma^+\in B_n(i^+)}p_{\sigma^+}=\sum_{j=1}^Nc^+_{j,i}.
\end{eqnarray*}
From (A1), we know that $P_1,P_2$ are both non-negative and irreducible. Thus, (\ref{s3}) is a consequence of Corollary 8.1.33 of \cite{HJ:13}.

Note that $\rho_1^n\to 0$ as $n\to\infty$. By  (\ref{s3}), for all large $n$, we have
\begin{equation}\label{guo6}
\sum_{\sigma\in B_n(i)}p_\sigma<\chi_i\chi_{i^+}^{-1}\underline{\chi}\sum_{\sigma^+\in B_n(i^+)}p_{\sigma^+},
\end{equation}
By Lemma \ref{condition1}, for every $\sigma\in B_n(i)$, we have
\begin{equation}  \label{guo3}
 \{\widetilde{\sigma}\in\Gamma(\sigma):\widetilde{\sigma}_n=i\}=\{\sigma\},\;\;
\{\widetilde{\sigma}\in\Gamma(\sigma):\widetilde{\sigma}_n=i^+\}\supset\{\sigma^+\}.
\end{equation}
Combining (\ref{guo6}) and (\ref{guo3}), we deduce
\begin{eqnarray*}
\sum_{\sigma\in B_n(i)}I_{1,\sigma}=\sum_{\sigma\in B_n(i)}\chi_{\sigma_1}p_\sigma<\chi_i\chi_{i^+}^{-1}
\sum_{\sigma^+\in B_n(i^+)}\chi_{\sigma^+_1}p_{\sigma^+}\leq \chi_i\chi_{i^+}^{-1}\sum_{\sigma\in B_n(i)}I_{2,\sigma}.
\end{eqnarray*}
It follows that there exists some $\sigma\in B_n(i)$ fulfilling (\ref{s4}).
\end{proof}

\textbf{Proof of Theorem \ref{mthm1} (2)}

We assume that (A1)-(A3) hold. By Lemma \ref{condition3}, for every $i\in\Psi$, there exists some $\sigma\in\mathcal{S}^*(i)$
such that $\chi_{i^+}I_{1,\sigma}<\chi_iI_{\sigma,2}$. Thus, by Theorem \ref{mthm1} (1), $\mu$ is reducible if and only if (\ref{sg1}) holds for every $(i,j)\in \mathcal{S}_2$.
Note that by (A1), we have that $p_{i^+,j}=0$ for every $(i,j)\in \mathcal{S}_2$. Therefore, $\mu$ is reducible if and only if
$p_{i,j}+p_{i,j^+}=p_{i^+,j^+}$ for every $(i,j)\in \mathcal{S}_2$.

Next, we construct some examples to illustrate our results and assumptions. Our first example shows that for a suitable transition matrix $P$ and some initial probability vector, the measure $\mu$ coincides with the in-homogeneous self-similar measure
that is studied in \cite{Zhu:16}.
\begin{example}\label{eg1}{\rm
Let $f_i,1\leq i\leq N$, be contractive similarity mappings on $\mathbb{R}^q$. Let $E$ be the self-similar set determined by $(f_i)_{i=1}^N$. Let $(q_i)_{i=0}^N$
and $(t_i)_{i=1}^N$ be two positive probability vectors. Let $\nu_0$ denote the self-similar measure associated with $(f_i)_{i=1}^N$ and $(t_i)_{i=1}^N$.
We define
\begin{eqnarray*}
P=\left(\begin{array}{cccccc}
q_1& \cdots &q_N &q_0t_1  &\cdots &q_0t_N\\
\vdots & &\vdots&\vdots& &\vdots\\
q_1 &\cdots &q_N&q_0t_1&\cdots&q_0t_N\\
0 &\cdots &0&t_1&\cdots&t_N\\\vdots & &\vdots&\vdots& &\vdots\\0 &\cdots &0&t_1&\cdots&t_N
\end{array}\right),\;\;\chi=\left(\begin{array}{c}
q_1\\
\vdots \\
q_N \\
q_0t_1 \\\vdots \\q_0t_N
\end{array}\right).
\end{eqnarray*}
For $1\leq i,j\leq N$, we define $T_{i,j}=T_{i,j^+}=T_{i^+,j^+}=f_i$.
Then the MW-fractal $K$ agrees with $E$.
By Lemma \ref{condition1} and (\ref{mu-meausre}), for $\sigma\in\Psi^n$, one easily gets
\begin{eqnarray*}
\mu(J_\sigma)=\prod_{h=1}^nq_{\sigma_h}+q_0\prod_{h=1}^nt_{\sigma_h}+q_0\sum_{h=1}^{n-1}\bigg(\prod_{l=1}^hq_{\sigma_h}\cdot\prod_{l=h+1}^nt_{\sigma_h}\bigg).
\end{eqnarray*}
Thus, $\mu$ agrees with the in-homogeneous self-similar measure in \cite{Zhu:16}. That is, the unique probability measures satisfying $\mu=q_0\nu_0+\sum_{i=1}^Nq_i\circ\mu\circ f_i^{-1}$.
As noted in \cite[Remark 1.4]{Zhu:16}, for fixed $r>0$ and suitably selected $P$, it indeed can happen that $s_{1,r}>s_{2,r}$, $s_{1,r}<s_{2,r}$, or $s_{1,r}=s_{2,r}$.
}\end{example}
Our second example shows that, if (A3) and (A4) are not satisfied, it can happen that
$ \inf_{\sigma\in\mathcal{S}^*}\frac{\mu(J_{\sigma^\flat})}{\mu(J_\sigma)}=0$, regardless of whether $P$ is reducible. In addition, it is possible that $\sigma,\tau\in\mathcal{S}^*$ and $(\sigma_{|\sigma|},\tau_1)\in\mathcal{S}_2$, but $\sigma\ast\tau\notin\mathcal{S}^*$. This might cause major difficulties in the estimation for the quantization errors.

\begin{example}\label{eg2}{\rm
For $N=3$, we assume that (\ref{overlapping mapps}) holds. We define
\begin{eqnarray*}
P(1)=\left(\begin{array}{cccccc}
\frac{1}{3}&\frac{1}{3} &0 &\frac{1}{3}&0 &0\\
0& \frac{1}{3} &\frac{1}{3} &0&0&\frac{1}{3}\\
\frac{1}{3} &0& \frac{2}{3}&0&0&0\\
0& 0 &0 &\frac{1}{3}&0&\frac{2}{3}\\
0& 0 &0 &0&\frac{1}{3}&\frac{2}{3}\\0&0&0&\frac{1}{3}& \frac{1}{3} &\frac{1}{3}
\end{array}\right),\;
P(2)=\left(\begin{array}{cccccc}
\frac{1}{3}&\frac{1}{3} &0 &\frac{1}{3}&0 &0\\
0& \frac{1}{3} &\frac{1}{3} &0&0&\frac{1}{3}\\
\frac{1}{3} &0& \frac{2}{3}&0&0&0\\
0& 0 &0 &\frac{1}{3}&0&\frac{2}{3}\\
\frac{1}{3}& 0 &0 &0&\frac{1}{3}&\frac{1}{3}\\0&0&0&\frac{1}{3}& \frac{1}{3} &\frac{1}{3}
\end{array}\right).
\end{eqnarray*}
Note that $P(1)$ is reducible, while $P(2)$ is irreducible, since
$1\rightarrow4\rightarrow6\rightarrow5\rightarrow1\rightarrow2\rightarrow3\rightarrow1$ forms a cycle in the graph $\mathcal{G}$. Let $\chi_1=\ldots=\chi_6=\frac{1}{6}$. For every $n\geq 1$, let $\sigma^{(n)}=(1,1,\ldots, 1)\in\mathcal{S}_n$. Either $P=P(1)$ or $P=P(2)$, one can see that $P_1$ and $P_2$ are both irreducible, and $p_{1,2}>0$ and $p_{1,5}=p_{4,2}=p_{4,5}=p_{4,1}=0$, implying that $(1,2)\in G_2$, but $(1,5),(4,2),(4,5),(4,1)\notin G_2$. Thus,
\[
\Gamma(\sigma^{(n)})=\{(1,\ldots,1),(1,4,4,4\ldots,4),(1,1,4,4,\ldots,4),\ldots,(4,4,\ldots,4)\},
\]
but $\Gamma(\sigma^{(n)}\ast 2)=\{\sigma^{(n)}\ast 2\}$. As one can see, $\frac{\mu(J_{\sigma^{(n)}\ast 2})}{\mu(J_{\sigma^{(n)}})}\to 0$ as $n\to\infty$. This also happens in case that $P=P(1)$, even if we add $(1,5)$ to $G_2$ by adjusting the first row of $P$.

Now let $P=P(2)$ and $\sigma=(1),\tau=(2,1)$. Since $(5,1)\in G_2$, we have, $(2,1)\in\mathcal{S}_2$. However, $(1,5),(2,1),(2,4),(4,5)\notin G_2$. This implies that $(1,2,1)\notin\mathcal{S}_3$.
}\end{example}

Our third example shows that even if (g2) holds, the measure $\mu$ may not be reducible, when both (b1) and
the condition in (g1) fail.
\begin{example}\label{eg3}{\rm
Let $N=3$. Let $P$ and $\chi=(\chi_i)_{i=1}^6$ be defined by
\begin{eqnarray*}
P=\left(\begin{array}{cccccc}
\frac{1}{6}&\frac{1}{3} &0 &\frac{1}{6}&\frac{1}{3} &0\\
0& \frac{1}{6} &\frac{1}{3} &0&\frac{1}{6} &\frac{1}{3}\\
\frac{1}{3} &0& \frac{1}{6}&\frac{1}{3} &0& \frac{1}{6}\\
\frac{1}{3}&\frac{1}{6}&0 &\frac{1}{3}&\frac{1}{6}&0\\
0& \frac{1}{3} &\frac{1}{6} &0&\frac{1}{3} &\frac{1}{6}\\\frac{1}{3} &0& \frac{1}{6}&\frac{1}{3} &0& \frac{1}{6}
\end{array}\right), \;\;\chi=\left(\begin{array}{c}
\frac{1}{6}\\
\frac{1}{9} \\
\frac{1}{6} \\
\frac{1}{6} \\\frac{2}{9} \\\frac{1}{6}
\end{array}\right).
\end{eqnarray*}
By Lemma \ref{T5}, one can see that
\[
\mathcal{S}_2=H_1^2=\{(1,1),(1,2),(2,2),(2,3),(3,1),(3,3)\}.
\]
For every $(l,i)\in\mathcal{S}_2$, we have, $p_{l,i}+p_{l^+,i}=p_{l,i^+}+p_{l^+,i^+}$. In addition,
\[
p_{2,3}+p_{2,6}\neq p_{5,3}+p_{5,6},\;\;\chi_2\neq\chi_5.
\]
Thus, both (g1) and (b1) fail, but the condition in (g2) is fulfilled.
Next, we show that $\mu$ is not reducible.
Let $\sigma=(1,2,3)$. The set $\Gamma(\sigma)$ is exactly given by
\begin{eqnarray*}
\{(1,2,3), (4,2,3),(1,5,3), (4,5,3),(1,2,6), (4,2,6), (1,5,6), (4,5,6)\}.
\end{eqnarray*}
By (\ref{mu-meausre}) and  (\ref{zsg5}), we easily get $\mu(J_\sigma)\neq\widetilde{\chi}_1\widetilde{p}_{1,2}\widetilde{p}_{2,3}$.
Hence, $\mu$ is not reducible.
}\end{example}

\section{Some estimates for the quantization error for $\mu$}

In this section, we assume that either (A1)-(A3) hold, or, (A2), (A4) and (A5) hold.
For $\sigma,\tau\in\mathcal{S}^*\cup\mathcal{S}_\infty$, We say that $\sigma$ is comparable with $\tau$ and write $\sigma\prec\tau$ if
$|\sigma|\leq|\tau|$ and $\sigma=\tau|_{|\sigma|}$.
 If we have neither $\sigma\prec\tau$ nor $\tau\prec\sigma$, then we say that $\sigma,\tau$ are incomparable. Write
 \[
 \underline{p}:=\min_{(i,j)\in G_2}p_{i,j};\;\overline{p}:=\max_{(i,j)\in G_2}p_{i,j};\;\;\underline{s}:=\min_{(i,j)\in \mathcal{S}_2}s_{i,j};\;\overline{s}:=\max_{(i,j)\in \mathcal{S}_2}s_{i,j}.
 \]
\begin{remark}\label{sg3}{\rm
Let $|A|$ denote the diameter of a set $A\subset\mathbb{R}^q$. Without loss of generality, we assume that $|J_i|=1$ for all $i\in\Psi$. Then using (\ref{overlapping mapps}), we have
\[
|J_\sigma|=s_\sigma,\;\sigma\in\mathcal{S}^*; \;\;s_{\widetilde{\sigma}}=s_\sigma,\;{\rm for\;\;all}\;\;\widetilde{\sigma}\in\Gamma(\sigma).
\]
We define $\mathcal{E}_r(\theta):=1$. For every $\sigma\in\mathcal{S}^*$, we define
\begin{eqnarray}\label{energy}
\mathcal{E}_r(\sigma):=\mu(J_\sigma) s_\sigma^r=
\sum_{\widetilde{\sigma}\in\Gamma(\sigma)}(\chi_{\widetilde{\sigma}_1}p_{\widetilde{\sigma}}s_{\widetilde{\sigma}}^r).
\end{eqnarray}
}
\end{remark}
Using the following lemmas, we present some basic facts about the cylinder sets and the measure $\mu$, so that Lemma 3 of \cite{KZ:16} is applicable.
\begin{lemma}\label{l3}
There exist some constants $c_{1,r},c_{2,r}\in(0,1)$, such that
\[
c_{1,r}\mathcal{E}_r({\sigma^\flat})\leq\mathcal{E}_r(\sigma)\leq c_{2,r}\mathcal{E}_r({\sigma^\flat}),\;\sigma\in\mathcal{S}^*.
\]
\end{lemma}
\begin{proof}
Note that $N\geq 2$. By (A2) and (A3), or, (A2) and (A4), we have
\begin{eqnarray*}\label{T7}
&&\overline{\zeta}:=\max_{1\leq i\leq N}(\chi_i+\chi_{i^+})<1;\;\;\min_{(i,j)\in \mathcal{S}_2}(p_{i^+,j}+p_{i^+,j^+})\geq\underline{p};
\\&&\overline{d}=\max_{(i,j)\in \mathcal{S}_2}\max\{(p_{i,j}+p_{i,j^+}),(p_{i^+,j}+p_{i^+,j^+})\}<1.
\end{eqnarray*}
For $n\geq 2$ and $\sigma\in\mathcal{S}_n$, we write $\sigma=\tau\ast j$. Then $\sigma^\flat=\tau$. By Lemma \ref{condition1} and (\ref{mu-meausre}), $\mu(J_i)=\chi_i+\chi_{i^+}$ for $i\in\Psi$; and for $\sigma\in\mathcal{S}_n$ with $n\geq 2$, we have
\begin{eqnarray*}\label{zg2}
\mu(J_\sigma)=I_{1,\tau} (p_{\sigma_{n-1},j}+p_{\sigma_{n-1},j^+})+I_{2,\tau}(p_{\sigma_{n-1}^+,j}+p_{\sigma_{n-1}^+,j^+}).
\end{eqnarray*}
The lemma follows by defining $c_{1,r}:=\min\{\underline{p},2\underline{\chi}\}\underline{s}^r$ and $c_{2,r}:=\max\{\overline{p},\overline{\zeta},\overline{d}\}\overline{s}^r$.
\end{proof}
\begin{lemma}\label{l3'}
Let $L\in\mathbb{N}$. There exist a constant $\delta>0$ and a number $D_L>0$ which is independent of $\sigma\in\mathcal{S}^*$  and $\alpha\subset\mathbb{R}^d$ such that
\begin{enumerate}
\item[(c1)] $d(J_\sigma,J_\tau)\geq \delta\max\{|J_\sigma|,|J_\tau|\}$ for incomparable words $\sigma,\tau\in\mathcal{S}^*$;

\item[(c2)]$\int_{J_\sigma}d(x,\alpha)^rd\mu(x)\geq D_L\mathcal{E}_r(\sigma)$ for $\alpha\subset\mathbb{R}^q$ of cardinality $L$ and $\sigma\in\mathcal{S}^*$.
\end{enumerate}
\end{lemma}
\begin{proof}
(c1) By (\ref{ssc}) and our assumption for $J_i, i\in\Psi$, for some constant $\delta>0$,
\begin{eqnarray}\label{ssc2}
&&d(J_i, J_j)\geq\delta\max\{|J_i|,|J_j|\},\;\;1\leq i\neq j\leq N;\\
&&d(J_{i,j}, J_{i,l})\geq\delta\max\{|J_{i,j}|,|J_{i,l}|\},\;\;j\neq l,\;\;(i,j),(i,l)\in \mathcal{S}_2.\nonumber
\end{eqnarray}
Let $\sigma,\tau\in\mathcal{S}^*$ be incomparable words. Let
$l=\min\{i\geq 1:\sigma_i\neq\tau_i\}$. We have that $J_\sigma\subset J_{\sigma|_l}$ and $J_\tau\subset J_{\tau|_l}$. If $l=1$, then (c2) follows by (\ref{ssc2}). For
 $l\geq 2$, we write $\sigma|_l=\rho\ast i$ and $\tau|_l=\rho\ast j$, for some $1\leq i\neq j\leq N$. Then
\begin{eqnarray*}\label{sgzhu2}
d(J_\sigma,J_\tau)\geq d(J_{\sigma|_l},J_{\tau|_l})\geq s_\rho\delta\max\{|J_{\rho_{l-1}\ast i}|,|J_{\rho_{l-1}\ast j}|\}\geq\delta\max\{|J_\sigma|,|J_\tau\}.
\end{eqnarray*}

(c2) Let $\sigma\in\mathcal{S}^*$ be given. By (A2), one can see that
\[
{\rm card}(\{\tau\in \mathcal{S}^*: \sigma\prec\tau, |\tau|=|\sigma|+h\})\geq 2^h.
\]
Thus, (c2) can be obtained from the proof of \cite[Lemma 4]{Zhu:08b}.
\end{proof}

For every $k\geq 1$ and $s\in(0,\infty)$, we define
\begin{eqnarray}\label{lambdakr}
&&\Lambda_{k,r}:=\{\sigma\in\mathcal{S}^*:\mathcal{E}_r(\sigma)<c_{1,r}^k\leq\mathcal{E}_r(\sigma^\flat)\};
\\&&\phi_{k,r}:={\rm card}(\Lambda_{k,r}),\;l_{1k}:=\min_{\sigma\in\Lambda_{k,r}}|\sigma|;\;\;l_{2k}:=\max_{\sigma\in\Lambda_{k,r}}|\sigma|;\nonumber\\
&&\underline{P}^s_r(\mu):=\liminf_{k\to\infty}\phi_{k,r}^{\frac{r}{s}}e_{\phi_{k,r},r}^r(\mu),\;
\overline{P}^s_r(\mu):=\limsup_{k\to\infty}\phi_{k,r}^{\frac{r}{s}}e_{\phi_{k,r},r}^r(\mu).\nonumber
\end{eqnarray}
\begin{remark}\label{g1}{\rm
The following facts will be useful in the proof of our main result:
\begin{enumerate}
\item[(d1)] $l_{1k},l_{2k}\asymp k$; This can be seen from the following facts:
\[
c_{1,r}^{l_{1k}}<c_{1,r}^k,\;{\rm and}\;c_{2,r}^{l_{2k}-1}\geq c_{1,r}^k.
\]
\item[(d2)] For $s\in(0,\infty)$, $\underline{Q}^s_r(\mu)>0$ if and only if $\underline{P}^s_r(\mu)>0$, and  $\overline{Q}^s_r(\mu)<\infty$ if and only if
    $\overline{P}^s_r(\mu)<\infty$ (cf. \cite[Lemma 2.4]{Zhu:13}).
\end{enumerate}
}\end{remark}

\begin{lemma}\label{l4a}
Assume that (A1)-(A3) hold, or, (A2), (A4), (A5) hold. Then
\[
e_{\phi_{k,r},r}^r(\mu)\asymp\sum\limits_{\sigma\in\Lambda_{k,r}}\mathcal{E}_r(\sigma).
\]
\end{lemma}
\begin{proof}
This follows from (\ref{lambdakr}), Lemma \ref{l3} and \cite[Lemma 3]{KZ:16}.
\end{proof}

For $s\in(0,\infty)$ and $k\geq 1$, we define $F^s_{k,r}(\mu):=\sum_{\sigma\in\Lambda_{k,r}}(\mathcal{E}_r(\sigma))^{\frac{s}{s+r}}$ and
\begin{eqnarray}\label{fkrs}
\underline{F}^s_r(\mu)=\liminf_{k\to\infty}F^s_{k,r}(\mu),\; \overline{F}^s_r(\mu)=\limsup_{k\to\infty}F^s_{k,r}(\mu).
\end{eqnarray}

Using techniques from \cite[Proposition 14.5, 14.11]{GL:00}, we are able to reduce the asymptotics of the quantization errors to those for the sequence $(F^s_{k,r}(\mu))_{k=1}^\infty$. That is,
\begin{lemma}\label{l4b} For every $s>0$, we have
\begin{enumerate}
\item[(1)] $\underline{Q}^s_r(\mu)>0$ if and only if $\underline{F}^s_r(\mu)>0$;
\item[(2)] $\overline{Q}^s_r(\mu)<\infty$ if and only if
$\overline{F}^s_r(\mu)<\infty$.
\end{enumerate}
\end{lemma}
\begin{proof}
Let $s>0$ be given. We only show (1), and (2) can be proved analogously.
Assume that $\underline{F}^s_r(\mu)=:\xi>0$. Then there exists some $k_1>0$ such that for every $k\geq k_1$, we have
$F^s_{k,r}(\mu)>\frac{\xi}{2}$. Using Lemma \ref{l4a} and H\"{o}lder's inequality with exponent less than one, we deduce
\begin{eqnarray*}
\phi_{k,r}^{\frac{r}{s}}e_{\phi_{k,r},r}^r(\mu)\gtrsim \phi_{k,r}^{\frac{r}{s}}\sum_{\sigma\in\Lambda_{k,r}}\mathcal{E}_r(\sigma)
\geq\phi_{k,r}^{\frac{r}{s}}\big(F^s_{k,r}(\mu)\big)^{\frac{s+r}{s}}\phi_{k,r}^{-\frac{r}{s}}
\geq\big(\frac{\xi}{2}\big)^{\frac{s+r}{s}}.
\end{eqnarray*}
It follows that $\underline{P}_r^s(\mu)>0$. This and Remark \ref{g1} (d2) yield that $\underline{Q}_r^s(\mu)>0$.

Now we assume that $\underline{F}^s_r(\mu)=0$. Then for every $\epsilon\in(0,1)$, there exists a subsequence $(k_i)_{i=1}^\infty$ of positive integers, such that $F^s_{k_i,r}(\mu)<\epsilon$
for every $i\geq 1$. By Lemma \ref{l3}, for every $\sigma\in\Lambda_{k_i,r}$, we have $\mathcal{E}_r(\sigma)\geq c_{1,r}^{k_i+1}$. It follows that
\[
\phi_{k_i,r}(c_{1,r})^{(k_i+1)s/(s+r)}\leq F^s_{k_i,r}(\mu)<\epsilon.
\]
Using this, (\ref{lambdakr}) and Lemma \ref{l4a}, we deduce
\begin{eqnarray*}\label{g2}
\phi_{k_i,r}^{\frac{r}{s}}e_{\phi_{k_i,r},r}^r(\mu)\lesssim
\phi_{k_i,r}^{\frac{r}{s}}\sum_{\sigma\in\Lambda_{k_i,r}}(\mathcal{E}_r(\sigma))^{\frac{s}{s+r}}c_{1,r}^{\frac{k_ir}{s+r}}<c_{1,r}^{-\frac{r}{s+r}}\epsilon.
\end{eqnarray*}
It follows that $\underline{P}_r^s(\mu)=0$. By Remark \ref{g1} (d2), we conclude that $\underline{Q}_r^s(\mu)=0$.
\end{proof}

\section{Proof of Theorem \ref{mthm}}

In this section, we assume that (A1)-(A3) hold.
 A subset $\Gamma$ of $\mathcal{S}^*$ is called a \emph{finite anti-chain} if $\Gamma$ is finite and words in $\Gamma$ are pairwise incomparable.
 A finite anti-chain $\Gamma$ is called \emph{maximal} if for every word $\tau\in \mathcal{S}_\infty$, there exists some $\sigma\in\Gamma$ such that $\sigma\prec\tau$.
 We define a finite (maximal) anti-chain in $G^*, H_i^*,i=1,2$, or $H_1^*(j), H_2^*(j^+)$ with $j\in\Psi$ analogously.
The following lemma provides us with a useful tool to estimate $F^s_{k,r}(\mu)$, which can be seen as a generalization of \cite[Lemma 3.1]{KZ:15}.
\begin{lemma}\label{spectrallem}
Let $\Gamma_1$ be an arbitrary finite maximal anti-chain  in $H_1^*$, or in $H_1^*(j)$ for some $j\in\Psi$ and $\Gamma_2$ a finite maximal anti-chain in $H_2^*$ or in $H_2^*(j^+)$. Let $l(\Gamma_i):=\min\limits_{\sigma\in\Gamma_i}|\sigma|$ and $L(\Gamma_i):=\max\limits_{\sigma\in\Gamma_i}|\sigma|, i=1,2$. Then for $s>0$, there exist positive numbers $c_5(s), c_6(s)$, which are independent of $\Gamma_i$, such that
\begin{eqnarray}\label{guo5}
\left\{\begin{array}{ll}c_5(s)\rho_i(s)^{l(\Gamma_i)}\leq\sum\limits_{\sigma\in\Gamma_i}(p_\sigma s_\sigma^r)^{\frac{s}{s+r}}\leq c_6(s)\rho_i(s)^{L(\Gamma_i)}&{\rm if}\;\;s\leq s_{i,r}\\
c_5(s)\rho_i(s)^{L(\Gamma_i)}\leq\sum\limits_{\sigma\in\Gamma_i}(p_\sigma s_\sigma^r)^{\frac{s}{s+r}}\leq c_6(s)\rho_i(s))^{l(\Gamma_i)}&{\rm if}\;\;s>s_{i,r}\end{array}\right..
\end{eqnarray}
\end{lemma}
\begin{proof}
It suffices to give the proof for $i=1$. Note that the spectral radius of $\rho_1(s)^{-1}A_1(\frac{s}{s+r})$ equals $1$.
Since $A_1(\frac{s}{s+r})$ is nonnegative and irreducible, by Perron-Frobenius theorem, there exists a unique positive normalized right eigenvector $(\xi_l)_{l=1}^N$ of
$\rho_1(s)^{-1}A_1(\frac{s}{s+r})$ with respect to $1$:
\[
\sum_{j=1}^N\rho_1(s)^{-1}(p_{l,j}s_{l,j}^r)^{\frac{s}{s+r}}\xi_j=\xi_l,\;1\leq l\leq N.
\]
For $k\geq 2$ and $\sigma\in H_1^k$, we define $\nu_1([\sigma]):=\rho_1(s)^{-k}(p_\sigma s_\sigma^r)^{\frac{s}{s+r}}\xi_{\sigma_k}$.
We have
\begin{eqnarray*}
\sum_{j=1}^N\nu_1([\sigma\ast j])&=&\rho_1(s)^{-k}(p_\sigma s_\sigma^r)^{\frac{s}{s+r}}\rho_1(s)^{-1}\sum_{j=1}^N(p_{\sigma_k,j}s_{\sigma_k,j}^r)^{\frac{s}{s+r}}\xi_j
\\&=&\rho_1(s)^{-k}(p_\sigma s_\sigma^r)^{\frac{s}{s+r}}\xi_{\sigma_k}=\nu_1([\sigma]).
\end{eqnarray*}
Thus, $\nu_1$ extends a measure on $H_1^{\infty}$. We distinguish the following two cases.
\begin{enumerate}
\item[\rm (1)] $\Gamma_1$ is a finite maximal anti-chain in $H_1^*$. We have
\begin{eqnarray}\label{g3}
\sum_{\sigma\in\Gamma_1}\nu_1(\Gamma_1)=\sum_{\sigma\in\Gamma_1}\rho_1(s)^{-|\sigma|}(p_\sigma s_\sigma^r)^{\frac{s}{s+r}}\xi_{\sigma_{|\sigma|}}
=\nu_1(H_1^\infty)=\rho_1(s)^{-1}.
\end{eqnarray}

\item[\rm (2)] $\Gamma_1$ is a finite maximal anti-chain in $H_1^*(j)$ for some $j\in\Psi$. We have
\begin{eqnarray}\label{g3'}
\sum_{\sigma\in\Gamma_1}\nu_1(\Gamma_1)=\sum_{\sigma\in\Gamma_1}\rho_1(s)^{-|\sigma|}(p_\sigma s_\sigma^r)^{\frac{s}{s+r}}\xi_{\sigma_{|\sigma|}}=\rho_1(s)^{-1}\xi_j.
\end{eqnarray}
\end{enumerate}
We define $\underline{\xi}:=\min\limits_{1\leq i\leq N}\xi_i$ and $\overline{\xi}:=\max\limits_{1\leq i\leq
N}\xi_i$. Then by (\ref{g3}), (\ref{g3'}), one can see that (\ref{guo5}) is fulfilled with
$c_5(s):=\overline{\xi}^{-1}\underline{\xi}\rho_1(s)^{-1}$ and $c_6(s):=\underline{\xi}^{-1}\rho_1(s)^{-1}$.
\end{proof}

For the proof of Theorem \ref{mthm}, we define
\begin{eqnarray}
&&\mathcal{A}_{k,r}:=\bigcup_{\sigma\in\Lambda_{k,r}}\bigg\{\sigma|_h\ast(\sigma_{h+1}^+,\ldots,\sigma_{|\sigma|}^+), 1\leq h\leq|\sigma|-1\bigg\};\label{Akr}
\\&& a_{k,r}(s):=\sum_{\sigma\in\Lambda_{k,r}}\sum_{h=1}^{|\sigma|-1}(p_{\sigma|_h\ast(\sigma_{h+1}^+,\ldots,\sigma_{|\sigma|}^+)}s_\sigma^r)^{\frac{s}{s+r}},\;k\geq1.
\nonumber
\end{eqnarray}

\textbf{Proof of Theorem \ref{mthm} (i)}

 By Lemma \ref{condition1}, $\{\sigma,\sigma^+\}\subset\Gamma(\sigma)$. Thus,
\begin{eqnarray*}
F^{s_r}_{k,r}(\mu)\geq\underline{\chi}^{\frac{s_r}{s_r+r}}\max\bigg\{\sum_{\sigma\in\Lambda_{k,r}}(p_\sigma s_\sigma^r)^{\frac{s_r}{s_r+r}},\sum_{\sigma\in\Lambda_{k,r}}(p_{\sigma^+}
s_{\sigma^+}^r)^{\frac{s_r}{s_r+r}}\bigg\}.
\end{eqnarray*}
Note that $\Lambda_{k,r}$ is a finite maximal anti-chain in $H_1^*$.
By (A3), $\{\sigma^+:\sigma\in\Lambda_{k,r}\}$ is a maximal anti-chain in $H_2^*$. By Lemma \ref{spectrallem}, one easily gets
$F^{s_r}_{k,r}(\mu)\gtrsim 1$. This and Lemma \ref{l4b} yield that $\underline{Q}_r^{s_r}(\mu)>0$ and $\underline{D}_r(\mu)\geq s_r$.
Next, we prove that $\overline{D}_r(\mu)\leq s_r$. For $1\leq h\leq l_{2k}-1$, we define
\[
B(\omega):=\{\tau^+\in H_2^*:\omega\ast\tau^+\in\mathcal{A}_{k,r}\},\;\omega\in H_1^h.
\]
For every $\omega\in H_1^h$, either $B(\omega)=\emptyset$ or $B(\omega)$ is an anti-chain in $H_2^*$. In fact, suppose that $\tau^+,\rho^+\in B(\omega)$ are distinct words with $\tau^+\prec\rho^+$; by (\ref{Akr}), we would have $\omega\ast\tau,\omega\ast\rho\in\Lambda_{k,r}$ which are comparable, a contradiction. We have
\begin{equation}\label{sanguo1}
\mathcal{A}_{k,r}\subset\bigcup_{h=1}^{l_{2k}-1}\bigcup_{\omega\in H_1^h}\{\omega\ast\tau^+:\tau^+\in B(\omega)\}.
\end{equation}
For $s>s_r$,  we have, $\rho_i(s)<1$ for $i=1,2$. By (\ref{sanguo1}) and Lemma \ref{spectrallem},
\begin{eqnarray}
F^s_{k,r}(\mu)&\leq&\sum_{\sigma\in\Lambda_{k,r}}(p_\sigma s_\sigma^r)^{\frac{s}{s+r}}+\sum_{\sigma\in\Lambda_{k,r}}(p_{\sigma^+} s_{\sigma^+}^r)^{\frac{s}{s+r}}+a_{k,r}(s)\label{zsg3}\\&\lesssim&2+\sum_{h=1}^{l_{2k}-1}\sum_{\omega\in H_1^h}\sum_{\tau^+\in B(\omega)}(p_\omega s_\omega^r)^{\frac{s}{s+r}}(p_{\tau^+} s_{\tau^+}^r)^{\frac{s}{s+r}}
\nonumber\\&\lesssim& 2+\frac{\rho_1(s)}{1-\rho_1(s)}.\label{guo4}
\end{eqnarray}
This and Lemma \ref{l4b}, yield that $\overline{Q}_r^s(\mu)<\infty$ and $\overline{D}_r(\mu)\leq s$. It follows that $\overline{D}_r(\mu)\leq s_r$. This completes the proof of Theorem \ref{mthm} (i).

Next, we are going to prove Theorem \ref{mthm} (ii). For $\tau^+\in H_2^*$, we define
\[
B(\tau^+):=\{\omega\in H_1^*:\omega\ast\tau^+\in\mathcal{A}_{k,r}\}=\{\omega\in H_1^*:\omega\ast\tau\in\Lambda_{k,r}\}.
\]
Clearly, $B(\tau^+)$ might be empty for some $\tau^+\in H_2^*$. Also, it can happen that two words $\omega,\rho\in B(\tau^+)$ are comparable while $\omega\ast\tau,\rho\ast\tau$ are
incomparable. In fact, this happens if for some $\upsilon\in\mathcal{S}^*$,
\[
\min\{\mathcal{E}_r(\omega\ast\tau^\flat),\mathcal{E}_r(\omega\ast\upsilon\ast\tau^\flat)\}\geq c_{1,r}^k>\max\{\mathcal{E}_r(\omega\ast\tau),\mathcal{E}_r(\omega\ast\upsilon\ast\tau)\}.
\]

The following Lemma \ref{l5}-\ref{l6} are devoted to the case that $s_{1,r}>s_{2,r}$. Our next lemma will enable us to estimate the difference $\big||\omega|-|\rho|\big|$ for any two comparable words $\omega, \rho\in B(\tau^+)$.
\begin{lemma}\label{l5}
Let $\omega\in H_{1,n},\tau\in H_{1,m},\upsilon\in H_{1,l}$ and $\omega\ast\tau, \omega\ast\upsilon\ast\tau\in \mathcal{S}^*$. Then
\[
\mathcal{E}_r(\omega\ast\upsilon\ast\tau)\leq (l+1)(\overline{p}\overline{s}^r)^{l+1}(\underline{p}\underline{s}^r)^{-1}\mathcal{E}_r(\omega\ast\tau).
\]
\end{lemma}
\begin{proof}
By Lemma \ref{condition1} and (\ref{mu-meausre}), with $\upsilon|_0:=\theta$, we have
\begin{eqnarray}
\mu(J_{\omega\ast\tau})&=&\chi_{\omega_1}p_{\omega\ast\tau}+\sum_{h=1}^{n-1}\chi_{\omega_1}p_{\omega|_h\ast(\omega_{h+1}^+,\ldots,\omega_n^+)\ast\tau^+)}
\label{e1}\\&&\;\;\;+\chi_{\omega_1}p_{\omega\ast\tau^+}\label{e2}\\&&\;\;\;+\sum_{h=1}^{m-1}\chi_{\omega_1}p_{\omega\ast\tau|_h\ast(\tau_{h+1}^+,\ldots,\tau_m^+)}
+\chi_{\omega^+_1}p_{\omega^+\ast\tau^+}\label{e3}.
\end{eqnarray}
Similarly, using Lemma \ref{condition1} and (\ref{mu-meausre}),  we have
\begin{eqnarray}
\mu(J_{\omega\ast\upsilon\ast\tau})&=&\chi_{\omega_1}p_{\omega\ast\upsilon\ast\tau}+
\sum_{h=1}^{n-1}\chi_{\omega_1}p_{\omega|_h\ast(\omega_{h+1
}^+,\ldots,\omega_n^+)\ast\upsilon^+\ast\tau^+)}\label{e4}\\&&\;\;\;
+\sum_{h=0}^l\chi_{\omega_1}p_{\omega\ast\upsilon|_h\ast(\upsilon_{h+1}^+,\ldots,\upsilon_l^+)\ast\tau^+}\label{e5}\\&&\;\;\;
+\sum_{h=1}^{m-1}\chi_{\omega_1}p_{\omega\ast\upsilon\ast\tau|_h\ast(\tau_{h+1}^+,\ldots,\tau_m^+)}+\chi_{\omega^+_1}p_{\omega^+\ast\upsilon^+\ast\tau^+}.\label{e6}
\end{eqnarray}
We denote the sum in (\ref{e1}), (\ref{e2}), (\ref{e3}), by $I_1, I_2, I_3$, and denote the sum in (\ref{e4}), (\ref{e5}), (\ref{e6}), by $I_4, I_5, I_6$. Then
\begin{eqnarray*}
I_1&=&\chi_{\omega_1}p_\omega\cdot p_{\omega_n,\tau_1}\cdot p_\tau
\\&&\;\;\;+\sum_{h=1}^{n-1}\chi_{\omega_1}p_{\omega|_h\ast(\omega_{h+1}^+,\ldots,\omega_n^+)}\cdot p_{\omega_n^+,\tau_1^+}\cdot  p_{\tau^+};\\
I_4&=&\chi_{\omega_1}p_{\omega}\cdot p_{\omega_n,\upsilon_1}\cdot p_\upsilon \cdot p_{\upsilon_1,\tau_1}\cdot p_\tau
\\&&\;\;\;+\sum_{h=1}^{n-1}\chi_{\omega_1}p_{\omega|_h\ast(\omega_{h+1}^+,\ldots,\omega_n^+)}\cdot p_{\omega_n^+,\upsilon_1^+}\cdot p_{\upsilon^+}\cdot p_{\upsilon_l^+,\tau_1^+}\cdot
p_{\tau^+}.
\end{eqnarray*}
We compare the preceding two equations and obtain
\begin{eqnarray}\label{compare1}
I_4\leq (\overline{p}^{l+1}\underline{p}^{-1})\cdot I_1.
\end{eqnarray}
In an similar manner, one can see that
\begin{eqnarray}\label{compare2}
I_6\leq (\overline{p}^{l+1}\underline{p}^{-1})\cdot I_3.
\end{eqnarray}
Next, we compare $I_5$ and $I_2$. We have
\begin{eqnarray}\label{compare3}
I_2&=&\chi_{\omega_1}p_{\omega\ast\tau^+}=\chi_{\omega_1}\cdot p_\omega\cdot p_{\omega_n,\tau_1^+}\cdot p_{\tau^+};\nonumber\\
I_5&=&\sum_{h=0}^l\chi_{\omega_1}p_\omega\cdot p_{\omega_n,\upsilon_1}\cdot p_{\upsilon|_h\ast(\upsilon_{h+1}^+,\ldots,\upsilon_l^+)}\cdot p_{\upsilon_l^+,\tau_1^+}\cdot
p_{\tau^+}\nonumber\\&\leq&(l+1)\overline{p}^{l+1}\chi_{\omega_1}\cdot p_\omega\cdot p_{\tau^+}\nonumber\\
&\leq&(l+1)\overline{p}^{l+1}\underline{p}^{-1} I_2.
\end{eqnarray}
Combining (\ref{compare1})-(\ref{compare3}), we deduce
\begin{eqnarray*}
\mu(J_{\omega\ast\upsilon\ast\tau})&=&I_4+I_5+I_6\leq (l+1)\overline{p}^{l+1}\underline{p}^{-1}(I_1+I_2+I_3)\\&\leq&(l+1)\overline{p}^{l+1}\underline{p}^{-1}\mu(J_{\omega\ast\tau}).
\end{eqnarray*}
Note that $s_{\omega\ast\upsilon\ast\tau}\leq \overline{s}^{l+1}\underline{s}^{-1}s_{\omega\ast\tau}$. We obtain
\[
\mathcal{E}_r(\omega\ast\upsilon\ast\tau)=\mu(J_{\omega\ast\upsilon\ast\tau})s_{\omega\ast\upsilon\ast\tau}^r
\leq(l+1)(\overline{p}\overline{s}^r)^{l+1}(\underline{p}\underline{s}^r)^{-1}\mathcal{E}_r(\omega\ast\tau).
\]
This completes the proof of the lemma.
\end{proof}

\begin{remark}\label{t2}{\rm
For $\tau^+\in H_2^*$, we define
\[
B^\flat(\tau^+):=\{\omega\in B(\tau^+):\rho\nprec\omega\;{\rm for\;\; every}\;\;\rho\in B(\tau^+)\setminus\{\omega\}\}.
\]
It is easy to see that $B^\flat(\tau^+)$ is an anti-chain in $H_1^*$.
}\end{remark}

For every $\omega\in B^\flat(\tau^+)$, we define $B_\omega(\tau^+):=\{\rho\in B(\tau^+):\omega\prec\rho\}$.
Next, we give an estimate for the size of $B_\omega(\tau^+)$ by using Lemma \ref{l5}. Let $k_2$ be the smallest integer such that
\[
(k_2+1)(\overline{p}\overline{s}^r)^{k_2+1}(\underline{p}\underline{s}^r)^{-1}<c_{1,r}.
\]
\begin{lemma}\label{l6}
 Let $C:=\sum_{h=0}^{k_2}N^h$. For every $\omega\in B^\flat(\tau^+)$, we have
 \begin{equation}\label{guo2}
 \sum_{\rho\in B_\omega(\tau^+)}(p_\rho s_\rho^r)^{\frac{s_{1,r}}{s_{1,r}+r}}\leq C(p_\omega s_\omega^r)^{\frac{s_{1,r}}{s_{1,r}+r}}.
 \end{equation}
\end{lemma}
\begin{proof}
Suppose that there exists some $\upsilon$ with $|\upsilon|>k_2$ such that $\omega\ast\upsilon\in B(\tau^+)$. Then by  Lemma \ref{l5}, we have
$\mathcal{E}_r(\omega\ast\upsilon\ast\tau)<c_{1,r}\mathcal{E}_r(\omega\ast\tau)$.
This contradicts (\ref{lambdakr}), because both $\omega\ast\tau$ and $\omega\ast\upsilon\ast\tau$ are elements of $\Lambda_{k,r}$. Thus,
\begin{eqnarray*}
B_\omega(\tau^+)\subset\{\rho\in H_1^*:\omega\prec\rho,|\rho|\leq |\omega|+k_2\};\;{\rm card}(B_\omega(\tau^+))\leq C.
\end{eqnarray*}
For every $\rho\in B_\omega(\tau^+)$, we have $p_\rho s_\rho^r\leq p_\omega s_\omega^r$. Thus, (\ref{guo2}) is fulfilled.
\end{proof}

  \textbf{Proof of Theorem \ref{mthm} (ii) and (iii)}

For Theorem \ref{mthm} (ii), we need to treat the following three cases.

\textbf{Case 1:} $s_{1,r}<s_{2,r}$.
In this case, we have $s_r=s_{2,r}$ and $\rho_1(s_{2,r})<1$. One can see that (\ref{guo4}) remains valid for $s=s_{2,r}>s_{1,r}$. Thus, we have
$\overline{Q}_r^{s_r}(\mu)<\infty$.

\textbf{Case 2:} $s_{1,r}>s_{2,r}$. By (\ref{Akr}) and (A3), $\rho\ast\tau^+\in\mathcal{A}_{k,r}$ if and only if $\rho\ast\tau\in\Lambda_{k,r}$. We have
\[
\mathcal{A}_{k,r}\subset\bigcup_{h=1}^{l_{2k}-1}\bigcup_{\tau^+\in H_2^h}\bigcup_{\omega\in B^\flat(\tau^+)}\{\rho\ast\tau^+:\rho\in B_\omega(\tau^+)\}.
\]
For $s=s_{1,r}>s_{2,r}$, by (\ref{zsg3}), Lemmas \ref{spectrallem}, \ref{l6} and  Remark \ref{t2},
we deduce
\begin{eqnarray*}\label{z3}
F^s_{k,r}(\mu)&\lesssim& 2+\sum_{h=1}^{l_{2k}-1}\sum_{\tau^+\in H_2^h}\sum_{\omega\in B^\flat(\tau^+)}\sum_{\rho\in B_\omega(\tau^+)}(p_{\rho\ast\tau^+}s_{\rho\ast\tau^+})^{\frac{s}{s+r}}\\
&\lesssim& 2+C\sum_{h=1}^{l_{2k}-1}\sum_{\tau^+\in H_2^h}(p_{\tau^+}s_{\tau^+}^r)^{\frac{s}{s+r}}\sum_{\omega\in B^\flat(\tau^+)}(p_\omega s_\omega^r)^{\frac{s}{s+r}}
\\&\lesssim&2+\frac{\rho_2(s)}{1-\rho_2(s)}.
\end{eqnarray*}
It follows that $\overline{F}_r(s)<\infty$. This and Lemma \ref{l4b} yield that $\overline{Q}_r^{s_r}(\mu)<\infty$.

\textbf{Case 3:} $s_{1,r}=s_{2,r}$. For every $1\leq h\leq l_{1k}-1$, $\Lambda_{k,r}(h):=\{\sigma|_h:\sigma\in\Lambda_{k,r}\}$ is a maximal anti-chain in $H_1^*$. Fix an arbitrary $\omega\in\Lambda_{k,r}(h)$. By (A2),
the set
$D_\omega:=\{\tau\in H_1^*:\omega\ast\tau\in\Lambda_{k,r}\}$ contains a maximal anti-chains in $H_1^*(j_1)$ for some $j_1\in\Psi$.
Hence, by (A3), the set $D^+_{\omega}:=\{\tau^+:\omega\ast\tau\in\Lambda_{k,r}\}$ contains some maximal anti-chain $\mathcal{A}(\omega)$ in $H_2^*(j_1^+)$, and $\{\omega\ast\tau^+:\tau^+\in
D^+_{\omega}\}\subset\mathcal{A}_{k,r}$. Thus,
 \begin{equation}\label{tem1}
 \mathcal{A}_{k,r}\supset\bigcup_{h=1}^{l_{1k}-1}\bigcup_{\omega\in\Lambda_{k,r}(h)}\{\omega\ast\tau^+:\tau^+\in\mathcal{A}(\omega)\}.
 \end{equation}
By Lemma \ref{condition1} and H\"{o}lder's inequality with exponent less than one, we have
 \begin{eqnarray*}
 F^{s_r}_{k,r}(\mu)&\geq&\sum_{\sigma\in\Lambda_{k,r}}\bigg(\sum_{\widetilde{\sigma}\in\Gamma(\sigma)\setminus\{\sigma,\sigma^+\}}
 \chi_{\widetilde{\sigma}_1}p_{\widetilde{\sigma}}s_\sigma^r\bigg)^{\frac{s_r}{s_r+r}}\\
 &\geq&\sum_{\sigma\in\Lambda_{k,r}}\sum_{\widetilde{\sigma}\in\Gamma(\sigma)\setminus\{\sigma,\sigma^+\}}
 \big(\chi_{\widetilde{\sigma}_1}p_{\widetilde{\sigma}}s_\sigma^r\big)^{\frac{s_r}{s_r+r}}|\widetilde{\sigma}|^{-\frac{r}{s_r+r}}
\\&\geq&(\underline{\chi}l_{2k})^{-\frac{r}{s_r+r}}\sum_{\hat{\sigma}\in\mathcal{A}_{k,r}}(p_{\hat{\sigma}} s_{\hat{\sigma}}^r)^{\frac{s_r}{s_r+r}}.
 \end{eqnarray*}
 Using this, (\ref{tem1}), Lemma \ref{spectrallem} and Remark \ref{g1} (d1), we deduce
 \begin{eqnarray*}
 F^{s_r}_{k,r}(\mu)&\gtrsim& l_{2k}^{-\frac{r}{s_r+r}}\sum_{h=1}^{l_{1k}-1}\sum_{\omega\in\Lambda_{k,r}(h)}(p_\omega s_\omega^r)^{\frac{s_r}{s_r+r}} \sum_{\tau^+\in\mathcal{A}(\omega)}
 (p_{\tau^+} s_{\tau^+}^r)^{\frac{s_r}{s_r+r}}
 \\&\gtrsim& l_{2k}^{-\frac{r}{s_r+r}}l_{1k}\asymp k^{\frac{s_r}{s_r+r}}.
 \end{eqnarray*}
 Thus, by Lemma \ref{l4b}, we conclude that $\underline{Q}_r^{s_r}(\mu)=\infty$.

(iii) Let $G(i^+):=\{j^+: (i^+,j^+)\in G_2\}$. By (A2), for $i\in\Psi$, we have ${\rm card}(G(i^+))\geq 2$.
Let $t_{i,r}, i\in\Psi$, be implicitly defined by
\[
\sum_{j^+\in G(i^+)}(p_{i^+,j^+}s_{i^+,j^+}^r)^{\frac{t_{i,r}}{t_{i,r}+r}}=1.
\]
By (A1), we have $P_4=\textbf{0}$. Thus, for every $i\in\Psi$, $(p_{i^+,j^+})_{j^+\in G(i^+)}$ is a probability vector. By Theorem 14.14 of \cite{GL:00},
$t_{i,r}=D_r(\lambda_{i^+})$ for the self-similar measure $\lambda_{i^+}$ associated with $(p_{i^+,j^+})_{j^+\in G(i^+)}$ and some IFS $(f_{i^+,j^+})_{j^+\in G(i^+)}$
with similarity ratios $(s_{i^+,j^+})_{j^+\in G(i^+)}$. Thus, for every $r>0$, we apply \cite[Theorem 11.6]{GL:00} and obtain
$t_{i,r}\geq\frac{\log \overline{p}}{\log \underline{s}}=:\kappa>0$.
Let $\zeta_r:=\min\limits_{1\leq i\leq N}t_{i,r}$. By Theorem 8.1.22 of \cite{HJ:13},
\[
\rho_2(s)\geq\min_{1\leq i\leq N}\sum_{j=1}^N(p_{i^+,j^+}s_{i^+,j^+}^r)^{\frac{s}{s+r}}=\min_{1\leq i\leq N}\sum_{j^+\in G(i^+)}(p_{i^+,j^+}s_{i^+,j^+}^r)^{\frac{s}{s+r}}.
\]
It follows that $s_{2,r}\geq \zeta_r\geq\kappa$. For every $i\in\Psi$, we have
\begin{eqnarray*}
\limsup_{r\to 0}\sum_{j=1}^N(p_{i,j}s_{i,j}^r)^{\frac{s_{2,r}}{s_{2,r}+r}}
\leq\sum_{j=1}^N\limsup_{r\to 0}p_{i,j}^{\frac{\kappa}{\kappa+r}}(s_{i,j}^r)^{\frac{\kappa}{\kappa+r}}=\sum_{j=1}^Np_{i,j}.
\end{eqnarray*}
By (A2) and (A3), we know that $\sum_{j=1}^Np_{i,j}<1$ for every $i\in\Psi$. Thus, there exists some $r_0>0$ such that for every $r\in(0,r_0)$, we have
$\rho_1(s_{2,r})<1$. It follows that $s_{2,r}>s_{1,r}$. This and Theorem \ref{mthm} (ii) yield Theorem \ref{mthm} (iii).

\section{Proof of Theorem \ref{mthm2}}
In this section, we always assume that (A2), (A4) and (A5) hold. For every $\sigma\in\mathcal{S}^*$, let $\mathcal{E}_r(\sigma)$ be as defined in (\ref{energy}). We have

\begin{lemma}\label{t6}
Let $\sigma\ast\tau\in\mathcal{S^*}$. For $s>0$, we have
\[
(\underline{p}\underline{s}^r\overline{\chi}^{-1})^s(\mathcal{E}_r(\sigma))^s(\mathcal{E}_r(\tau))^s\leq
\mathcal{E}_r(\sigma\ast\tau)^s\leq(\overline{p}\overline{s}^r\underline{\chi}^{-1})^s(\mathcal{E}_r(\sigma))^s
(\mathcal{E}_r(\tau))^s.
\]
\end{lemma}
\begin{proof}
By Lemma \ref{T5}, for every $\sigma\in\mathcal{S}^*$, we have
\begin{equation*}\label{guo7}
\Gamma(\sigma\ast\tau)=\mathcal{T}(\sigma\ast\tau)=\{\widetilde{\sigma}\ast\widetilde{\tau}:\widetilde{\sigma}\in\mathcal{T}(\sigma),
\widetilde{\tau}\in\mathcal{T}(\tau)\}.
\end{equation*}
By (\ref{overlapping mapps}), for every $\rho\in\mathcal{S}^*$ and $\widetilde{\rho}\in\Gamma(\rho)$, we have $s_{\widetilde{\rho}}=s_\rho$.  Hence,
\begin{eqnarray*}
\big(\mathcal{E}_r(\sigma\ast\tau)\big)^s&=&\bigg(\sum_{\widetilde{\sigma}\in\mathcal{T}(\sigma)}\sum_{\widetilde{\tau}\in\mathcal{T}(\tau)}
\chi_{\widetilde{\sigma}_1}(p_{\widetilde{\sigma}}s_{\widetilde{\sigma}}^r)\cdot p_{\widetilde{\sigma}_n,\widetilde{\tau}_1}s_{\widetilde{\sigma}_n,\widetilde{\tau}_1}^r\cdot
(p_{\widetilde{\tau}}s_{\widetilde{\tau}}^r)\bigg)^s
\\&\geq&(\underline{p}\underline{s}^r\overline{\chi}^{-1})^s\bigg(\sum_{\widetilde{\sigma}\in\mathcal{T}(\sigma)}
(\chi_{\widetilde{\sigma}_1}p_{\widetilde{\sigma}}s_{\widetilde{\sigma}}^r)\sum_{\widetilde{\tau}\in\mathcal{T}(\tau)}
(\chi_{\widetilde{\tau}_1}p_{\widetilde{\tau}}s_{\widetilde{\tau}}^r)\bigg)^s\\
&=&(\underline{p}\underline{s}^r\overline{\chi}^{-1})^s(\mathcal{E}_r(\sigma))^s
(\mathcal{E}_r(\tau))^s.
\end{eqnarray*}
The remaining part of the lemma can be obtained similarly.
\end{proof}

For every $s>0$ and $n\geq 1$, we define
\begin{equation}\label{tns}
T_n(s):=\sum_{\sigma\in\mathcal{S}_n}\big(\mathcal{E}_r(\sigma)\big)^s=\sum_{\sigma\in\mathcal{S}_n}\bigg(\sum_{\widetilde{\sigma}\in\mathcal{T}(\sigma)}\chi_{\widetilde{\rho}_1}
p_{\widetilde{\sigma}}s_{\widetilde{\sigma}}^r\bigg)^s.
\end{equation}
Next, we are going to show that $(T_n(s))_{n=1}^\infty$ is quasi-multiplicative up to constant factors. For $h,l\geq 1$ and $\sigma\in\mathcal{S}_h$, we write
\[
\Lambda(\sigma,l):=\{\rho\in\mathcal{S}_l: \sigma\ast\rho\in\mathcal{S}_{h+l}\},\;\mathcal{S}_{l,i}:=\{\tau\in \mathcal{S}_l:\tau_1=i\},\;i\in\Psi.
\]
For $s>0, l\geq 1$ and $i\in\Psi$, we define
\[
T_{l,i}(s):=\sum_{\omega\in\mathcal{S}_{l,i}}\big(\mathcal{E}_r(\omega)\big)^s,\;h(s):=(Nc_{2,r}^s)^{-N}(\overline{p}\overline{s}^r\underline{\chi}^{-1})^{-s}(c_{1,r})^{Ns}(\underline{p}\underline{s}^r\overline{\chi}^{-1})^s.
\]
\begin{lemma}\label{temp4}
 Let $s>0$ be given. There exist positive numbers $g_1(s), g_2(s)$ such that for every pair $n,l\in\mathbb{N}$, we have
\begin{equation}\label{tem6}
g_1(s)T_n(s)T_l(s)\leq T_{n+l}(s)\leq g_2(s)T_n(s)T_l(s).
\end{equation}
\end{lemma}
\begin{proof}
We first show the following claim:

\textbf{Claim}: for $1\leq i\neq j\leq N$, we have $T_{l,i}(s)\geq h(s)T_{l,j}(s)$. By (A5), $P_1$ is irreducible, so the sub-graph $\mathcal{G}_1$ of $\mathcal{G}$ with vertex set $\Psi$ is strongly connected. Hence, there exists a word $\gamma$ with $|\gamma|<N$, such that $i\ast\gamma\ast j\in H_1^*=\mathcal{S}^*$.
Note that $\{i\ast\gamma\ast\tau:\tau\in \mathcal{S}_{l,j}\}\subset\mathcal{S}_{l+|\gamma|+1,i}$. Using Lemmas \ref{l3} and \ref{t6}, we deduce
\begin{eqnarray*}\label{zhu1}
T_{l+|\gamma|+1,i}(s)\left\{\begin{array}{ll}
\geq\sum\limits_{\tau\in\mathcal{S}_{l,j}}
\big(\mathcal{E}_r(i\ast\gamma\ast\tau)\big)^s \geq(c_{1,r})^{Ns}(\underline{p}\underline{s}^r\overline{\chi}^{-1})^s
T_{l,j}(s)\\
=\sum\limits_{\omega\in\mathcal{S}_{l,i}}\sum\limits_{\upsilon\in\Lambda(\omega,|\gamma|+1)}(\mathcal{E}_r({\omega\ast\upsilon}))^s
\leq (Nc_{2,r}^s)^{N}(\overline{p}\overline{s}^r\underline{\chi}^{-1})^sT_{l,i}(s)\end{array}\right..
\end{eqnarray*}
This completes the proof of the claim.

Now let $g_1(s):=N^{-1}h(s)(\underline{p}\underline{s}^r\overline{\chi}^{-1})^{s}$ and $g_2(s):=(\overline{p}\overline{s}^r\underline{\chi}^{-1})^s$. For every $\sigma\in\mathcal{S}_n$ and $l\geq 1$, we have, $\Lambda(\sigma,l)\subset\mathcal{S}_l$. Let $j_1\in\Psi$ such that $\mathcal{S}_{l,j_1}\subset\Lambda(\sigma,l)$. Using Lemma \ref{t6} and the claim, we deduce
\begin{eqnarray}\label{zg3}
\sum\limits_{\omega\in\Lambda(\sigma,l)}(\mathcal{E}_r(\sigma\ast\omega))^s\left\{\begin{array}{ll}
\leq\sum\limits_{\tau\in\mathcal{S}_l}(\mathcal{E}_r(\sigma\ast\tau))^s
\leq g_2(s)(\mathcal{E}_r(\sigma))^sT_l(s)\\
\geq\sum\limits_{\omega\in\mathcal{S}_{l,j_1}}\big(\mathcal{E}_r(\sigma\ast\omega)\big)^s
\geq g_1(s)\big(\mathcal{E}_r(\sigma)\big)^sT_l(s)\end{array}\right..
\end{eqnarray}
Since $\mathcal{S}_{n+l}=\bigcup\limits_{\sigma\in\mathcal{S}_n}\{\sigma\ast\omega:\omega\in\Lambda(\sigma,l)\}$, (\ref{tem6}) follows from (\ref{zg3}).
\end{proof}
\begin{remark}{\rm
Let $\sigma\in\mathcal{S}_n, j\in\Psi$ with $\sigma_n\ast j\in\mathcal{S}_2$. If (A4) is not satisfied, then it can happen that $\mathcal{S}_{l,j}\nsubseteq\Lambda(\sigma,l)$. This can be seen from Example \ref{eg2}. Let $P=P(2)$ and $\sigma=(1,1), l=2, j_1=2, \tau=(2,1)$. As we have noted, $\sigma\in\mathcal{S}_2$ and $\sigma_2\ast j_1=(1,2)\in\mathcal{S}_2$. Also, we have $\tau\in\mathcal{S}_2$ and $\tau\in\mathcal{S}_{2,2}$, but $\sigma\ast\tau\notin\mathcal{S}_4$. This means that $\tau\notin\Lambda(\sigma,2)$.
}\end{remark}

\begin{lemma}\label{zhu3}
For $s>0$, the limit $\lim\limits_{n\to\infty} \frac{1}{n}\log T_n(s)=:\Phi(s)$ exists. Moreover,
\begin{enumerate}
\item [(f1)] $\Phi$ is continuous and strictly decreasing;
there exists a unique $s_0\in(0,1)$ and a unique $t_r>0$ such that $\Phi(s_0)=0$ and $\Phi(\frac{t_r}{t_r+r})=0$;
\item[(f2)]for $b:=g_2(s_0)/g_1(s_0)$ and for every pair $m,n\in\mathbb{N}$, we have
\begin{equation*}\label{temp5}
 b^{-1}T_n(s_0)\leq T_m(s_0) \leq bT_n(s_0) .
\end{equation*}
\end{enumerate}
\end{lemma}
\begin{proof}
The lemma can be proved by using (A2), Lemmas \ref{l3}, \ref{temp4} and \cite[Corollary 1.2]{Fal:97} along the line of \cite[Lemma 5.2]{Fal:97}.
\end{proof}

Using some ideas in \cite[Theorem 5.1]{Fal:97}, we are now able to obtain an auxiliary measure. As we will see, this measure is closely connected with the quantization errors for $\mu$.
\begin{lemma}\label{guo1}
There exists a probability measure $\lambda$ supported on $K$, such that
\[
\lambda(J_\sigma)\asymp(\mathcal{E}_r(\sigma))^{\frac{t_r}{t_r+r}},\;\;\sigma\in\mathcal{S}^*.
\]
\end{lemma}
\begin{proof}
For $m\geq 1$ and $\sigma\in\mathcal{S}_m$, let $x_\sigma$ be an arbitrary point of $J_\sigma\cap K$ and denote by $\delta_\sigma$ the Dirac measure at the point $x_\sigma$.
For every $m\geq 1$, we define
$\lambda_m:=\frac{1}{T_m(s_0)}\sum_{\sigma\in\mathcal{S}_m}(\mathcal{E}_r(\sigma))^{\frac{t_r}{t_r+r}}\delta_\sigma$.
Then $(\lambda_m)_{m=1}^\infty$ is a sequence of probability measures.
By  \cite[Theorem 1.23]{PM:95}, there exist a sub-sequence $(\lambda_{m_k})_{k=1}^\infty$ and a measure $\lambda$ such that $\lambda_{m_k}\to\lambda$ (weak convergence)
as $k\to\infty$. One can see that $\lambda$ is a probability measure supported on $K$. Now let $n\geq 1$ and $\sigma\in\mathcal{S}_n$ be given. For every $m>n$,
\begin{eqnarray*}
\lambda_m(J_\sigma)=\sum_{\rho\in\Lambda(\sigma,m-n)}\lambda_m(\sigma\ast\rho)=\sum_{\rho\in\Lambda(\sigma,m-n)}\frac{1}{T_m(s_0)}
(\mathcal{E}_r({\sigma\ast\rho}))^{\frac{t_r}{t_r+r}}.
\end{eqnarray*}
Using (\ref{zg3}) and Lemma \ref{zhu3} (f2), one can easily obtain
\begin{eqnarray*}
 b^{-1}g_1(s_0) (\mathcal{E}_r({\sigma}))^{\frac{t_r}{t_r+r}}\leq\lambda_m(J_\sigma)
\leq bg_2(s_0) (\mathcal{E}_r({\sigma}))^{\frac{t_r}{t_r+r}}.
\end{eqnarray*}
This implies the assertion of the lemma .
\end{proof}

For the proof for Theorem \ref{mthm2} (ii), we define
\[
\widehat{p}_{i,j}=\left\{\begin{array}{ll}
p_{i,j}+p_{i,j^+},&{\rm in \;Case\; (g1)}\\
p_{i,j}+p_{i^+,j},&{\rm in \;Case\; (g2)}\end{array}\right.,\;(i,j)\in\Psi^2.
\]
For every $i\in\Psi$, $(\widehat{p}_{i,j})_{j=1}^N$ is a probability vector. In fact, we have
\[
\sum_{j=1}^N\widehat{p}_{i,j}=\left\{\begin{array}{ll}
\sum_{j=1}^N(p_{i,j}+p_{i,j^+})=\sum_{j=1}^{2N}p_{i,j}=1,&{\rm Case\; (g1)}\\
\frac{1}{2}\sum_{j=1}^N(p_{i,j}+p_{i^+,j}+p_{i,j^+}+p_{i^+,j^+})=1&{\rm Case\; (g2)}\end{array}\right..
\]
Let $B(s):=((\widehat{p}_{i,j}s_{i,j}^r)^s)_{i,j=1}^N$. We denote by $\xi(s)$ the spectral radius of $B(s)$. By (A2), there exists a unique $a_r>0$ such
that $\xi(\frac{a_r}{a_r+r})=1$. We have
\begin{lemma}\label{add}
Assume that (g1) or (g2) holds. We have $t_r=a_r$.
\end{lemma}
\begin{proof}
For every $\sigma\in\mathcal{S}_n$ and $1\leq h\leq n-1$, we define
\[
E_{\sigma_h,\sigma_{h+1}}:=\left(\begin{array}{cc}
p_{\sigma_h,\sigma_{h+1}}s_{\sigma_h,\sigma_{h+1}}^r&p_{\sigma_h,\sigma_{h+1}^+}s_{\sigma_h,\sigma_{h+1}}^r\\
p_{\sigma_h^+,\sigma_{h+1}}s_{\sigma_h,\sigma_{h+1}}^r& p_{\sigma_h^+,\sigma_{h+1}^+}s_{\sigma_h,\sigma_{h+1}}^r\\
\end{array}\right).
\]
Let $V:=(1\; 1)$ and $U=:V^T$. Let $\|x\|_1$ denote the $l_1$-norm for $x\in\mathbb{R}^2$. We have
\begin{equation}\label{zg5}
\sum_{\widetilde{\sigma}\in\mathcal{T}(\sigma)}
p_{\widetilde{\sigma}}s_{\widetilde{\sigma}}^r
=\big\|\prod_{h=1}^{n-1}E_{\sigma_h,\sigma_{h+1}}U\big\|_1=\big\|V\prod_{h=1}^{n-1}E_{\sigma_h,\sigma_{h+1}}\big\|_1.
\end{equation}
 Using (\ref{zg5}) and the conditions in (g1) and (g2), we deduce
\[
(\mathcal{E}_r(\sigma)^{\frac{s}{s+r}}\asymp\left\{\begin{array}{ll}
\|\prod_{h=1}^{n-1}E_{\sigma_h,\sigma_{h+1}}U\|_1^{\frac{s}{s+r}}=2^{\frac{s}{s+r}}(\widehat{p}_\sigma s_\sigma^r)^{\frac{s}{s+r}} &{\rm Case\; (g1)}\\
\|V\prod_{h=1}^{n-1}E_{\sigma_h,\sigma_{h+1}}\|_1^{\frac{s}{s+r}}
=2^{\frac{s}{s+r}}(\widehat{p}_\sigma s_\sigma^r)^{\frac{s}{s+r}}&{\rm  Case\; (g2)}\end{array}\right..
\]
Therefore, $t_r$ satisfies $\Upsilon(t_r):=\lim\limits_{n\to\infty}\frac{1}{n}\log\sum\limits_{\sigma\in\mathcal{S}_n}(\widehat{p}_\sigma s_\sigma^r)^{\frac{t_r}{t_r+r}}=0$.
As we did for $\Phi(s)$, one can see that the solution of $\Upsilon(s)=0$ is unique.
By applying Lemma \ref{spectrallem}, we know that
$\sum_{\sigma\in\mathcal{S}_n}(\widehat{p}_\sigma s_\sigma^r)^{\frac{a_r}{a_r+r}}\asymp 1$, implying that $a_r=t_r$.
\end{proof}

\textbf{Proof of Theorem \ref{mthm2}}

 (i) For every $k\geq 1$, by Lemma \ref{guo1}. We have
 \begin{eqnarray*}
 F^{t_r}_{k,r}(\mu)=\sum_{\sigma\in\Lambda_{k,r}}(\mathcal{E}_r(\sigma))^{\frac{t_r}{t_r+r}}\asymp\sum_{\sigma\in\Lambda_{k,r}}\lambda(J_\sigma)=1.
 \end{eqnarray*}
 This and Lemma \ref{l4b} yield that $0<\underline{Q}_r^{t_r}(\mu)\leq\overline{Q}_r^{t_r}(\mu)<\infty$ and $D_r(\mu)=t_r$.
 Note that $\bigcup_{\sigma\in\mathcal{S}_n}\Gamma(\sigma)=G_n$ is a finite maximal anti-chain in $G^*$ and that $P,A(\frac{s_r}{s_r+r})$ are irreducible.
By Lemma \ref{spectrallem}, we deduce
 \begin{eqnarray*}
T_n\big(\frac{s_r}{s_r+r}\big)=\sum_{\sigma\in\mathcal{S}_n}
\bigg(\sum_{\widetilde{\sigma}\in\Gamma(\sigma)}(\chi_{\widetilde{\sigma}_1}p_{\widetilde{\sigma}}s_{\widetilde{\sigma}}^r)\bigg)^{\frac{s_r}{s_r+r}}<\sum_{\widetilde{\sigma}\in G_n}(p_{\widetilde{\sigma}}s_{\widetilde{\sigma}}^r)^{\frac{s_r}{s_r+r}}\asymp 1.
 \end{eqnarray*}
 Hence, $\Phi(\frac{s_r}{s_r+r})\leq 0$ and $s_r\geq t_r$. This completes the proof of Theorem \ref{mthm2} (i).

(ii) Assume that (g1) or (g2) holds. By Lemma \ref{add}, we have $D_r(\mu)=t_r=a_r$. Next, we show that $s_r>a_r$.
Since $P_1$ is irreducible, so is the matrix $B(\frac{a_r}{a_r+r})$. There exists a positive right eigenvector $v=(v_1,\ldots,v_N)^T$ of $B(\frac{a_r}{a_r+r})$
in case (g1) and a positive left eigenvector $w=(w_1,\ldots,w_N)$ of $B(\frac{a_r}{a_r+r})$ in Case (g2), with respect to eigenvalue $1$:
\[
\left\{\begin{array}{ll}
\sum_{j=1}^N\big(\widehat{p}_{i,j}s_{i,j}^r\big)^{\frac{a_r}{a_r+r}}v_j=v_i&{\rm Case\; (g1)}\\
\sum_{j=1}^N\big(\widehat{p}_{j,i}s_{j,i}^r\big)^{\frac{a_r}{a_r+r}}w_j=w_i&{\rm  Case\; (g2)}\end{array}\right.,\;1\leq i\leq N.
\]
Let $\widetilde{v}=(v_1,\ldots v_N,v_1,\ldots, v_N)^T$ and $\widetilde{w}=(w_1,\ldots w_N,w_1,\ldots, w_N)$. Then $\widetilde{v},\widetilde{w}$ are positive vectors.
For every  $1\leq i\leq 2N$, let $R_i$ denote the $i$th row of the matrix $A(\frac{a_r}{a_r+r})$ and $C_i$ its $i$th column. By (\ref{overlapping mapps}), for every $(i,j)\in \mathcal{S}_2$, we have $s_{i,j}=s_{i,j^+}=s_{i^+,j}=s_{i^+,j^+}$. Hence,
\[
\left\{\begin{array}{ll}
R_i\widetilde{v}=\sum_{j=1}^N\big(p_{i,j}s_{i,j}^r\big)^{\frac{a_r}{a_r+r}}v_j
+\sum_{j=1}^N\big(p_{i,j^+}s_{i,j^+}^r\big)^{\frac{a_r}{a_r+r}}v_j>v_i&{\rm Case\; (g1)}\\
\widetilde{w}C_i=\sum_{j=1}^N\big(p_{j,i}s_{j,i}^r\big)^{\frac{a_r}{a_r+r}}w_j
+\sum_{j=1}^N\big(p_{j^+,i}s_{j^+,i}^r\big)^{\frac{a_r}{a_r+r}}w_j>w_i&{\rm Case\; (g2)}\end{array}\right..
\]
Using (g1) and (g2), one can also see that $(R_{i^+})\widetilde{v}>v_i$ and $\widetilde{w}(C_{i^+})>w_i$, for every $i\in\Psi$. It follows that
$A(\frac{a_r}{a_r+r})\widetilde{v}>\widetilde{v}$ in Case (g1) and $\widetilde{w}A(\frac{a_r}{a_r+r})>\widetilde{w}$ in Case (g2).
Thus, by applying \cite[Corollary 8.1.29] {HJ:13}, we obtain that $\rho(a_r)>1$ and $s_r>a_r$ . This completes the proof of Theorem \ref{mthm2} (ii).

\begin{remark}\label{equivalence}{\rm
In Case (g1), $\mu$ agrees with the Markov-type measure $\widehat{\mu}$ associated with $(\widehat{p}_{i,j})_{i,j=1}^N$ and $\widehat{\chi}=(\chi_i+\chi_{i^+})_{i=1}^N$. As Example \ref{eg3} shows, in Case (g2), it may happen that $\mu$ is not of Markov-type; however, $\mu$ is equivalent to $\widehat{\mu}$, in the sense that $\mu(A)\asymp\widehat{\mu}(A)$ for all Borel sets $A$. In fact, we have
\[
\mu(J_\sigma)=\sum_{\widetilde{\sigma}\in\Gamma(\sigma)}\chi_{\widetilde{\sigma}_1}
p_{\widetilde{\sigma}}\asymp\sum_{\widetilde{\sigma}\in\Gamma(\sigma)}p_{\widetilde{\sigma}}=2\widehat{p}_\sigma\asymp\widehat{\mu}(J_\sigma).
\]

Let $\widetilde{p}_{i,j},i,j\in\Psi$, be as defined in (\ref{zsg5}). Let $\widetilde{\mu}$ be the Markov-type measure associated with $(\widetilde{p}_{i,j})_{i,j=1}^N$ and $\widetilde{\chi}=(\chi_i+\chi_{i^+})_{i=1}^N$. As the following example shows, when $\mu$ is not reducible, $\mu$ and $\widetilde{\mu}$ are, in general, not equivalent.
}\end{remark}
\begin{example}\label{eg5}{\rm
Let $N=3$. Let $P$ be defined by
\begin{eqnarray*}
P=\left(\begin{array}{cccccc}
\frac{1}{6}&\frac{1}{6} &0 &\frac{1}{3}&\frac{1}{3} &0\\
0& \frac{1}{6} &\frac{1}{6} &0&\frac{1}{3} &\frac{1}{3}\\
\frac{1}{6} &0& \frac{1}{6}&\frac{1}{3} &0& \frac{1}{3}\\
\frac{1}{4}&\frac{1}{4}&0 &\frac{1}{8}&\frac{3}{8}&0\\
0& \frac{1}{4} &\frac{1}{4} &0&\frac{1}{4} &\frac{1}{4}\\\frac{1}{4} &0& \frac{1}{4}&\frac{1}{4} &0& \frac{1}{4}
\end{array}\right).
\end{eqnarray*}
Then (A2) and (A5) hold, but  for $i=1$, neither (a) nor (b) in Theorem \ref{mthm1} (1) holds, provided that $\chi_1=\chi_4$. In fact, $(1,1),(3,1)\in\mathcal{S}^*(1)$, and $(1,1),(1,4),(4,1),(4,4)\in G_2$, but
\[
p_{1,1}+p_{1,4}\neq p_{4,1}+p_{4,4};\;\;p_{1,1}+p_{4,1}\neq p_{1,4}+p_{4,4}.
\]
Thus, the measure $\mu$ is not reducible, when $\chi_1=\chi_4$. Next, we further show that $\mu$ is not even equivalent to the Markov-type measure $\widetilde{\mu}$, regardless of the choice of $\chi$. We consider
\[
\sigma^{(n)}=(1,1,\ldots, 1)\in\mathcal{S}_n,\;\;\tau^{(n)}=(1,2,3,1,2,3,\ldots, 1,2,3,1)\in\mathcal{S}_{3n+1}.
\]
Let $R_1$ denote the spectral radius of the following matrix:
\[
M_{11}:=\left(\begin{array}{cc}
p_{1,1}&p_{1,4}\\
p_{4,1}& p_{4,4}\\
\end{array}\right)=\left(\begin{array}{cc}
\frac{1}{6}&\frac{1}{3}\\
\frac{1}{4}& \frac{1}{8}\\
\end{array}\right).
\]
We have, $R_1=\frac{1}{96}(14+\sqrt{772})>0.43525$. Let
$\alpha_{i,j}^{(n)}$ denote the $(i,j)$-entry of the matrix $M_{11}^{n-1}$.
By Corollary 8.1.33 of \cite{HJ:13}, $\alpha_{1,1}^{(n)}+\alpha_{1,2}^{(n)}\asymp R^{n-1}$ and $\alpha_{2,1}^{(n)}+\alpha_{2,2}^{(n)}\asymp R_1^{n-1}$.
  We deduce
\begin{eqnarray}\label{eg22}
 \mu(J_{\sigma^{(n)}})\asymp
 \sum_{\widetilde{\sigma}\in\Gamma(\sigma^{(n)})}p_{\widetilde{\sigma}}=\sum_{i,j=1,2}\alpha_{i,j}^{(n)}\asymp R_1^{n-1},\;\widetilde{\mu}(\sigma^{(n)})\asymp\widetilde{p}_{1,1}^{n-1}.
\end{eqnarray}
Let $R_2$ denote the spectral radius of the following matrix:
\[
M_{1231}:=\left(\begin{array}{cc}
\frac{1}{6}&\frac{1}{3}\\
\frac{1}{4}& \frac{3}{8}\\
\end{array}\right)\cdot\left(\begin{array}{cc}
\frac{1}{6}&\frac{1}{3}\\
\frac{1}{4}& \frac{1}{4}\\
\end{array}\right)\cdot\left(\begin{array}{cc}
\frac{1}{6}&\frac{1}{3}\\
\frac{1}{4}& \frac{1}{4}\\
\end{array}\right).
\]
We have $R_2>0.1428$. By (\ref{zsg5}), we have $\widetilde{p}_{2,3}=\widetilde{p}_{3,1}=2^{-1}$. Thus,
\begin{equation}\label{eg23}
\mu(\tau^{(n)})=\|M_{1231}^n\|_1\asymp R_2^n,\;\widetilde{\mu}(\tau^{(n)})\asymp (\widetilde{p}_{1,2}\widetilde{p}_{2,3}\widetilde{p}_{3,1})^n\asymp 4^{-n}\widetilde{p}_{1,2}^n.
\end{equation}
Now we assume that $\mu$ and $\widetilde{\mu}$ are equivalent. Then by (\ref{eg22}) and (\ref{eg23}), we obtain
$\widetilde{p}_{1,1}=R_1$ and $\widetilde{p}_{1,2}=4R_2$. Setting $\zeta:=\chi_1/
(\chi_1+\chi_4)$, we have
\[
\frac{1}{2}\zeta+\frac{3}{8}(1-\zeta)=R_1;\;\;\frac{1}{2}\zeta+\frac{5}{8}(1-\zeta)=4R_2.
\]
However, we have $8(R_1-0.375)>0.482$, but $8(0.625-4R_2)<0.431$, a contradiction. Therefore, $\mu, \widetilde{\mu}$ are not equivalent, regardless of the choice of $\chi$.

}\end{example}
%\noindent{\bf Acknowledgements} This work was supported by National Natural Science Foundation of China (Grant No.
%11571144).

\end{document}